\documentclass{amsart}

\usepackage{graphicx}
\usepackage{amssymb,mathrsfs}
\usepackage{subfigure}

\usepackage{enumerate}        

\usepackage[all,dvips]{xy}    
\UseComputerModernTips        

\newcommand {\bng}{B_n\Gamma}

\newcommand {\cng}{\mathcal{C}^n\Gamma}
\newcommand {\ucng}{U\mathcal{C}^n\Gamma}
\newcommand {\ucn}[1]{U\mathcal{C}^n {#1}}
\newcommand {\dng}{\mathcal{D}^n\Gamma}
\newcommand {\udng}{U\mathcal{D}^n\Gamma}
\newcommand {\udn}[1]{U\mathcal{D}^n {#1}}
\newcommand {\bn}[1]{B_n #1}
\newcommand {\uc}[2]{U\mathcal{C}^{#1} #2}
\newcommand {\ud}[2]{U\mathcal{D}^{#1} #2}

\newtheorem{theorem}{Theorem}[section]
\newtheorem{lemma}[theorem]{Lemma}

\newtheorem{proposition}[theorem]{Proposition}
\newtheorem{corollary}[theorem]{Corollary}
\newtheorem{conjecture}[theorem]{Conjecture}
\theoremstyle{definition}
\newtheorem{definition}[theorem]{Definition}
\newtheorem{example}[theorem]{Example}

\begin{document}

\title[Cohomology of Tree Braid Groups]{Presentations for the Cohomology Rings of tree braid groups}
\author[D.Farley]{Daniel Farley}
      \address{Department of Mathematics \\
               Miami University of Ohio \\
               Oxford, OH 45056 \\}
      \email{farleyds@muohio.edu}

\begin{abstract}
If $\Gamma$ is a finite graph and $n$ is a natural number, then 
$\ucng$, the \emph{unlabelled configuration space of $n$ points on 
$\Gamma$}, is the space of all $n$-element subsets of $\Gamma$.  The
fundamental group of $\ucng$ is the \emph{$n$-strand braid group of
$\Gamma$}, denoted $\bng$.

We build on earlier work to compute presentations of the integral
cohomology rings $H^{\ast}( \uc{n}{T}; \mathbb{Z} )$, where $T$ is any tree
and $n$ is arbitrary.  The results suggest that 
$H^{\ast}( \uc{n}{T}; \mathbb{Z})$ is the exterior face ring of a
simplicial complex.
\end{abstract}



\maketitle

\section{Introduction} \label{intro}

Fix a finite graph $\Gamma$ and a natural number $n$.  The 
\emph{labelled configuration space of $n$ points on $\Gamma$},
denoted $\cng$, is the space
$$ \left( \prod_{i=1}^{n} \Gamma \right) - \Delta,$$
where $\Delta$ is the collection of ordered $n$-tuples 
$(x_1 , \ldots, x_n )$ such that $x_i = x_j $ for some $i \neq j$.
The \emph{unlabelled configuration space}, denoted $\ucng$, is the quotient
of $\cng$ by the action of the symmetric group $S_n$, where $S_n$ permutes
the factors.  
The fundamental groups of these spaces are, respectively,
the \emph{pure braid group} $P\bng$ and the \emph{braid group} $\bng$
of $n$ strands on the graph $\Gamma$.

This paper continues the study of the braid groups $\bng$ begun in \cite{FS1}.
Here we completely compute the cohomology ring 
$H^{\ast}(\ucn{T}; \mathbb{Z} )$, where $T$ is a tree.  The spaces
$\ucng$ are aspherical for any $n$ and any $\Gamma$ \cite{Ghr}, so, 
as an immediate result, we also compute 
the cohomology ring $H^{\ast}(\bn{T}; \mathbb{Z})$ of any 
\emph{tree braid group}, i.e., any braid group $\bn{T}$, where $T$ is a tree. 

In \cite{FS2}, 
Lucas Sabalka and I used partial information about the cohomology
rings $H^{\ast}( \bn{T}; \mathbb{Z}/2 \mathbb{Z})$ of tree braid groups
in order to characterize those groups $\bn{T}$ that are right-angled
Artin groups.  The crucial idea in computing the ring structure is to
introduce a ``cup product complex" $\widehat{\udn{T}}$, which is obtained
from a certain cubical complex $\udn{T}$ 
(described in Subsection \ref{basics}) 
by identifying all opposite faces of all cubes $C$ in $\udn{T}$.
The main argument of \cite{FS2}: i) showed that $\widehat{\udn{T}}$ is a
subcomplex of a high-dimensional torus; ii) showed that the map
$q^{\ast} : H^{\ast}( \widehat{\udn{T}}; \mathbb{Z}/2\mathbb{Z})  
\rightarrow H^{\ast}( \udn{T}; \mathbb{Z}/2\mathbb{Z})$ induced by 
the quotient map $q: \udn{T} \rightarrow \widehat{\udn{T}}$ is surjective,
and iii) explicitly described $q^{\ast}$.  As a result, it was possible
to compute $a \cup b$ for many elements 
$H^{\ast}( \udn{T}; \mathbb{Z}/ 2\mathbb{Z})$.

Here we achieve two improvements.  First, we work with integral coefficients.  
Second, we give a complete set of
generators for the ideal $\mathrm{Ker} \, q^{\ast}$, where $q^{\ast}$ is the map on
cohomology induced by the quotient map 
$q: \udn{T} \rightarrow \widehat{\udn{T}}$.  The argument, in outline,
goes like this.  The quotient map $q: \udn{T} \rightarrow \widehat{\udn{T}}$
induces a map $q^{\ast}: C^{\ast}( \widehat{\udn{T}}; \mathbb{Z})
\rightarrow C^{\ast}( \udn{T}; \mathbb{Z})$ on cochains.  All boundary maps
in $C_{\ast}( \widehat{\udn{T}}; \mathbb{Z})$ are $0$, so cohomology classes
in $H^{\ast}(\widehat{\udn{T}}; \mathbb{Z})$ are essentially the same as
cochains in $C^{\ast}( \widehat{\udn{T}}; \mathbb{Z})$.  The images of the
standard basis for $C^{\ast}(\widehat{\udn{T}}; \mathbb{Z})$ under the map $q^{\ast}$
are thus
cellular cocycles, which we call \emph{standard cocycles}.  The standard
cocycles can be easily visualized using what we call \emph{cloud pictures}.
A certain type of standard cocycle is called \emph{critical}; these form
a complete set of representatives for a basis of 
$H^{\ast}( \udn{T}; \mathbb{Z})$.  We describe an algorithm which takes a
standard cocycle $\phi$ as input and gives a sum $\Sigma$ of critical cocycles
as output, where $\Sigma$ and $\phi$ are cohomologous.  Moreover,
the difference $\Sigma - \phi$ is a sum of coboundaries from a certain
family, which we describe quite explicitly using another type of picture.
This easily gives the presentation promised by the title.

The paper concludes with the computation of several examples, which suggest
that the cohomology rings $H^{\ast}(\udn{T}; \mathbb{Z})$ are all
exterior face rings (see Subsection \ref{ghrist} for a definition).  

A large part of this paper was written while I was a guest 
at the Max Planck Institute for
Mathematics.  It is a pleasure to thank the Institute for
its hospitality and the agreeable working conditions during my stay.
 
\section{General Background on Configuration Spaces of Graphs} \label{gb} 

\subsection{Configuration spaces of graphs and graph braid groups}
\label{basics}

Let $\Gamma$ be a finite graph, and fix a natural number $n$. 
The \emph{labelled
configuration space} of $n$ points on $\Gamma$, denoted $\cng$, is the space
  $$\left(\prod^{n} \Gamma \right) - \Delta,$$
where $\Delta$ is the set of all points $(x_1, \dots, x_n) \in \prod^n
\Gamma$ such that $x_i = x_j$ for some $i \neq j$. The \emph{unlabelled
configuration space} of $n$ points on $\Gamma$, 
denoted $\ucng$, is the quotient of the
labelled configuration space by the action of the symmetric group $S_n$,
where the action permutes the factors. 
The \emph{$n$-strand braid group} of $\Gamma$,
denoted $\bng$, is the fundamental group of the
unlabelled configuration space of $n$ points on $\Gamma$. 


It will be convenient to use a space somewhat smaller than 
$\ucng$ in our applications.
Let $\Delta'$ denote the union of those open cells of $\prod^n \Gamma$
whose closures intersect $\Delta$. Let $\dng$ denote the
space $\left( \prod^n \Gamma \right) - \Delta'$. 
Note that $\dng$ inherits a CW complex
structure from the Cartesian product, and that an open cell in $\dng$ has the
form $c_1 \times \dots \times c_n$, where each $c_i$ is either a vertex
or the interior of an edge, and the closures of the $c_i$ are mutually
disjoint. Let $\udng$ denote the quotient of $\dng$ by the action of the
symmetric group $S_n$, which permutes the coordinates. An open cell in
$\udng$ has the form $\{c_1, \dots, c_n\}$, where each $c_i$ is either
a vertex or the interior of an edge, and the closures of any two of 
the $c_i$ are 
disjoint.  The set notation is used to indicate that order does not
matter.

In many circumstances, 
$\cng$ (respectively, $\ucng$) is homotopy equivalent to
$\dng$ (respectively, $\udng$).  Specifically:

\begin{theorem} \cite{A} \label{thm:Abrams}
For any $n>1$ and any graph $\Gamma$ with at least $n$ vertices, the
labelled (unlabelled) configuration space of $n$ points on $\Gamma$ strong
deformation retracts onto $\dng$ ($\udng$) if
\begin{enumerate}
\item each path between distinct vertices of degree not equal to $2$
passes through at least $n-1$ edges; and
\item each path from a vertex to itself which is not null-homotopic in
$\Gamma$ passes through at least $n+1$ edges.
\end{enumerate}
\end{theorem}

A graph $\Gamma$ satisfying the conditions of this theorem for a given $n$
is called \emph{sufficiently subdivided} for this $n$.  It is clear that
every graph is homeomorphic to a graph that is sufficiently 
subdivided for $n$, no matter
what $n$ may be.

We will sometimes call a vertex of degree $3$ or more \emph{essential}. 

Throughout the rest of the paper, we work exclusively with the space
$\udng$, where $\Gamma$ is sufficiently subdivided for $n$.  Also, from now
on, ``edge'' and ``cell'' will mean  ``closed edge" and ``closed cell",
respectively.

\subsection{An ordering on vertices of $\Gamma$ and a classification of
cells in $\udng$} \label{conventions} 

\begin{figure} [!h]
\begin{center}
\includegraphics{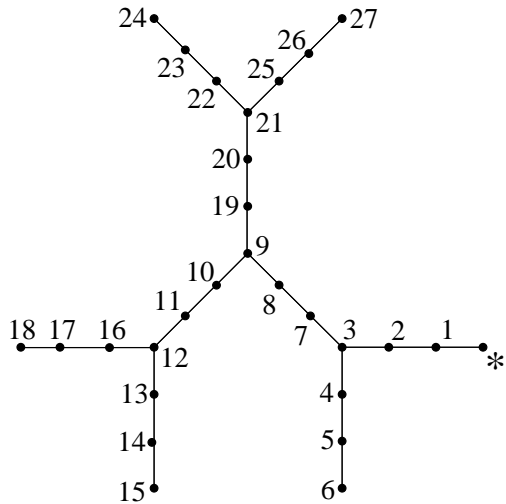}
\caption{This figure depicts the correct numbering of the given tree for the
given embedding.  Note that the tree is sufficiently subdivided for $n=4$ 
in the sense
of Subsection \ref{basics}.}
\label{number}
\end{center}
\end{figure}

Choose a maximal tree $T$ in $\Gamma$.   
Pick a vertex $\ast$ of valence $1$ in $T$ to be
the root of $T$. Choose an embedding of the tree $T$ into the plane.  We
define an order on the vertices of $T$ (and, thus, on vertices of
$\Gamma$) as follows. Begin at the basepoint $\ast$ and walk along the
tree, following the leftmost branch at any given intersection, and
consecutively number the vertices in the order in which they are first
encountered.  (When you reach a vertex of degree one, turn around.)  The
vertex adjacent to $\ast$ is assigned the number $1$.  
See Figure \ref{number}.  
Note that this
numbering depends only on the choice of $\ast$ and the embedding of the
tree.  

Let $\iota(e)$ and $\tau(e)$ denote the endpoints of a given edge
$e$ of $\Gamma$.  We orient each edge to go
from $\iota(e)$ to $\tau(e)$, and so that $\iota(e) < \tau(e)$.  (Thus, if
$e \subseteq T$, the geodesic segment $[ \tau(e), \ast ]$ in $T$ must pass
through $\iota(e)$.  Note that this convention for orientation is the 
opposite of that in \cite{Me} and \cite{FS1}.)

We use the order on the vertices to 
classify the cells $c$ of $\udng$ as collapsible, critical, or redundant.
These terms come from discrete Morse theory.  Our arguments will use
this theory only indirectly, so we omit a general introduction,
which can be found in \cite{For}.  
Strictly speaking, our cells are collapsible,
critical, or redundant relative to a discrete gradient vector field $W$
on $\ud{n}{T}$, which depends only on the choices of embedding of $T$ 
in the plane and basepoint $\ast$.  A thorough account of the discrete
gradient vector field $W$ on $\ud{n}{T}$ can be found in \cite{FS1}.

We require several definitions.
Edges outside of $T$ are called
\emph{deleted edges}.  
If $c = \{ c_1, \ldots, c_{n-1}, v \}$ and $e$ is the unique
edge in $T$ such that $\tau(e) = v$, then $v$ is \emph{blocked} with
respect to $c$ if $v= \ast$ or $\{ c_1, \ldots, c_{n-1}, e \}$ is not a
cell of $UD^{n} \Gamma$, i.e., if $c_{i} \cap e \neq \emptyset$ for some
$i \in \{ 1, \ldots, n-1 \}$.  Otherwise, $v$ is \emph{unblocked}.  If $c =
\{ c_1, \ldots, c_{n-1}, e \}$, the edge $e$ is \emph{not 
order-respecting} (with respect to $c$) if 
\begin{enumerate}
  \item there is a vertex $v$ in $c$ such that 
  \begin{enumerate}
    \item $v$ is adjacent to $\iota(e)$, and
    \item $\iota(e) < v < \tau(e)$, or
  \end{enumerate}    
  \item $e$ is a deleted edge.
\end{enumerate}
Otherwise, the edge $e$ is \emph{order-respecting} with respect to $c$.

It will often be useful to have another definition.  If $v$ is a
vertex in the tree $T$, we say that two vertices $v_1$ and $v_2$ lie in
the same \emph{direction} from $v$ if the geodesics $[ v , v_1 ] , [v ,
v_2 ] \subseteq T$ coincide in some neighborhood of $v$.  It follows that
there are
$n$ directions from a vertex of degree $n$ in $T$.  We number these
directions $0, 1 , 2, \ldots, n-1$, beginning with the direction
represented by $[v, \ast]$, which is numbered $0$, and proceeding in clockwise
order.  We will sometimes write $g( v_1 , v_2 )$ (where $v_1 \neq v_2$)  
to refer to the direction from $v_1$ to $v_2$.

Suppose that we are given a cell $c = \{ c_1, \ldots, c_n \}$ in $UD^{n}
\Gamma$.  Assign each cell $c_i$ in $c$ a number as follows.  A vertex of $c$ is
given the number from the above traversal of $T$.  An edge $e$ of $c$ is
given the number for $\tau(e)$.  Arrange the cells of $c$ in a sequence
$\mathcal{S}$, from the least- to the greatest-numbered.  (Here
the basepoint $\ast$ is numbered $0$.)  

\begin{definition}  
We classify each cell $c$ of $\ud{n}{\Gamma}$ 
as follows:
\begin{enumerate}
\item If an unblocked vertex occurs in $\mathcal{S}$ before all of the
order-respecting edges in $c$ (if any), then 
$c$ is redundant.
\item If an order-respecting edge occurs before any unblocked vertex, then
$c$ is collapsible.  
\item If all vertices of $c$ are blocked 
and no edge of $c$ is order-respecting, 
then $c$ is critical.
\end{enumerate}
(Of course, this definition is ultimately a theorem.  See Theorem 3.6 of
\cite{FS1}.) 
\end{definition}   

\begin{example} 
We give a short example to illustrate the previous definitions.  Consider
the  `Y'-graph, which is homeomorphic to the capital letter $Y$.  We consider
the discretized configuration space $\ud{3}{Y}$ for $Y$ sufficiently 
subdivided.  Figure \ref{ccr} depicts three different closed $1$-cells
in $\ud{3}{Y}$.  

\begin{figure} [!h]
\begin{center}
\includegraphics{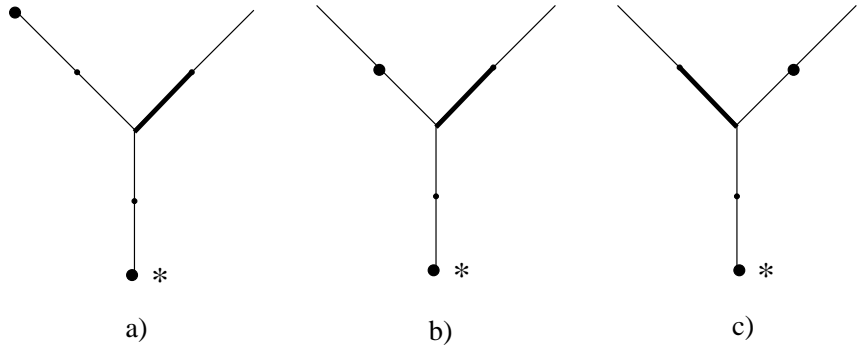}
\caption{These are three different $1$-cells in $\ud{3}{Y}$.  The cells are
redundant, critical, and collapsible (respectively).}
\label{ccr}
\end{center}
\end{figure}

In a), the vertices (large dots in the figure)
 are numbered $0$ and $4$, and the edge 
(represented by a thick line segment) 
is numbered $5$.  The vertex $0$ is blocked.  The vertex $4$ is
unblocked, so the cell in a) is redundant.  Note that the 
edge is order-respecting.

In b), the vertices are numbered $0$ and $3$, and the edge
is again numbered $5$.  The vertex $0$ is again blocked, and so
is the vertex $3$.  The edge $e$, numbered $5$, is not order-respecting,
since the vertex numbered $3$ lies between the initial and terminal
vertices of $e$ (numbered $2$ and $5$, respectively), and vertex $3$ is
adjacent to vertex $2$.  It follows from this that the $1$-cell
in b) is critical.

In c), the vertices are numbered $0$ and $5$, and the edge is numbered $3$.
The vertex $0$ is blocked.  The edge $e$ is order-respecting, since
$\iota(e)$ is numbered $2$, and $\tau(e)$ is numbered $3$, and thus
it is impossible to satisfy condition (1b) in the above definition
of order-respecting edges.  
(This particular edge must be order-respecting
in any cell $c$ containing it.)  It follows that the cell in c) is
collapsible.  Note that the vertex numbered $5$ is blocked.  \qed 
\end{example}

\subsection{On the homology of tree braid groups}

Discrete Morse theory (which is implicit in the previous subsection) allows
one to compute the homology groups of a complex $X$ 
from a simplified chain complex
$(M_{i}(X), \tilde{\partial}_{i})$, called
the \emph{Morse complex}.  A basis for the $i$th chain group
$M_{i}(X)$ is given by the collection of critical $i$-cells (that is, 
critical relative
to a discrete gradient vector field $W$ -- see \cite{For}).  The boundary maps
$\tilde{\partial}_{i}$ have natural definitions, but we won't need them  
here.

In \cite{Me}, we gave a 
computation of the integral homology groups of the spaces
$\ud{n}{T}$ in the form of the following theorem.
 
\begin{theorem} \cite{Me}
Let $( M_{i}( \ud{n}{T}) , \tilde{\partial}_{i})$ be the Morse complex
relative to the discrete gradient vector field $W$ (as defined in 
\cite{FS1}).
Each map $\tilde{\partial}_{i}$ is $0$.  In particular, the $i$th integral
homology group of $\ud{n}{T}$ is free abelian of rank $c(i)$, where $c(i)$
is the number of critical $i$-cells in $\ud{n}{T}$.   
\qed
\end{theorem}

In particular, we can simply identify $H_{i}( \ud{n}{T} )$ with the
free abelian group generated by critical $i$-cells, and we do this freely
from now on.  In some cases, we will nevertheless
need to have a description of 
homology representatives in the standard cellular chain complex
$( C_{i}(\ud{n}{T}), \partial_{i})$.  The following Proposition will be
sufficient.

\begin{proposition} \cite{Me2}
Let $X$ be a finite regular
CW complex endowed with a discrete gradient vector field $W$.  The
isomorphism between $H_{i}(M(X))$ and $H_{i}(C(X))$ is induced by a
map $\phi_{i}: M_{i}(X) \rightarrow C_{i}(X)$ on the chain level, such
that $\phi_{i}(c) = c + (\mathrm{collapsible} \, \mathrm{cells})$, for any
critical $i$-cell $c$. 
\qed
\end{proposition}

Since we will be working with integral homology and cohomology, we will 
need a convention for orienting cells of $\udn{T}$.  Consider
a $k$-cell $c= \{ e_{1} , e_{2}, \ldots, e_{k}, v_{k+1}, \ldots, v_{n} \}$.
Suppose, without loss of generality, that $\iota(e_{1}) < \iota(e_{2}) <
\ldots < \iota(e_{k})$.  We identify $c$ with $[0,1]^{k}$ by a map
$f: [0,1]^{k} \rightarrow c$ which sends $( t_{1}, \ldots, t_{k} )$ to the point
$\{ e_{1}(t_{1}), e_{2}(t_{2}), \ldots, e_{k}(t_{k}), v_{k+1}, \ldots, v_{n} \}$.  (Note
that the set notation here represents an actual $n$-element subset of $T$ rather than
a collection of cells.)   
Here the map $e_{i}: [0,1] \rightarrow T$ is injective, and 
satisfies $e_{i}(0) = \iota(e_{i})$, $e_{i}(1) = \tau(e_{i})$, and $e_{i}([0,1]) = e_{i}$.
We give $[0,1]^{k}$ the standard orientation, and $c$ the orientation induced by the identification $f$.

\subsection{The cellular boundary maps in the cubical complexes $\udn{T}$}
Consider the \emph{standard $n$-cube} $[0,1]^{n}$.  The \emph{standard orientation}
for this standard $n$-cube is a certain choice of orthogonal frame, namely the ordered
$n$-tuple $( v_{1}, \ldots, v_{n})$, where $v_i$ is the $n$-dimensional vector which is
$0$ in every coordinate except for the $i$th one, where it is $1$.

The \emph{boundary map} $\partial : C_{n} \left( [ 0,1 ]^{n} \right) \rightarrow C_{n-1} \left( [0,1]^{n} \right)$
can be described as follows \cite{BG}:
$$ \partial \left( [0,1]^{n} \right) = \sum_{i=1}^{n} (-1)^{n} \left( A_{i} - B_{i} \right), $$
where $A_i$ is the $(n-1)$-dimensional cube $[0,1] \times [0,1] \times \ldots \times \{ 0 \} \times \ldots \times [0,1]$
in which the factor $\{ 0 \}$ occurs in the $i$th place.  The $(n-1)$-dimensional cube $B_i$ is the same, except that the
$i$th factor is $\{ 1 \}$ instead of $\{ 0 \}$.  All of the cubes $A_i$, $B_i$ are given the standard orientations
$( v_1 , v_2 , \ldots , v_{i-1} , v_{i+1}, \ldots, v_{n} )$.

In a general cubical complex $X$ each cube $c$ must be given an orientation by specifying an identification $p: [0,1]^{k} \rightarrow
c$.  The boundary $\partial (c)$ in the cellular chain complex for $X$ may then be computed from the above formula, using suitable
identifications, as determined by $p$.  With the conventions adopted at the end of Subsection 
\ref{conventions}, we get the following lemma.

\begin{lemma} \label{basform}
Let $c = \{ e_1 , e_2 , \ldots, e_{k}, v_{k+1}, \ldots, v_{n} \}$ be a $k$-cell in $\udn{T}$,
where $\iota(e_{1}) < \iota(e_{2}) < \iota(e_{3}) < \ldots < \iota(e_{k})$, and $T$ is a tree.
Suppose that all cells of $\udn{T}$ are given their standard orientations (as in \ref{conventions}).
$$ \partial (c) = \sum_{i=1}^{k} (-1)^{i} \left( c_{\iota(e_{i})} - c_{\tau(e_{i})} \right),$$
where $c_{\iota(e_{i})}$ is the $(k-1)$-cube obtained from $c$ by replacing the edge $e_i$ with
its initial vertex $\iota(e_{i})$.  The $(k-1)$-cube $c_{\tau(e_{i})}$ is defined analogously, and both
are given their standard orientations.
\end{lemma}

\begin{proof}
This is a matter of checking definitions and the consistency of the orientations.
\end{proof}

The following lemma, an immediate consequence of the previous one, is true because of our very careful choice
of orientations on the cubes of $\udn{T}$.    

\begin{lemma} \label{miracle}  
Suppose that 
$c = \{ e_{1}, \ldots, e_{i-1}, e_{i}, e_{i+1}, \ldots, e_{k}, v_{k+1}, \ldots, v_{n} \}$,
$c' = \{ e_{1}, \ldots, e_{i-1}, e'_{i}, e_{i+1}, \ldots, e_{k}, v_{k+1}, \ldots, v_{n} \}$,
$\iota(e_{1})  < \iota(e_{2}) < \iota(e_{3}) < \ldots < \iota(e_{k})$, and $\iota(e_{i}) = \iota(e'_{i})$.
The cell  $c_{\iota(e_{i})} = c'_{\iota(e'_{i})}$ occurs in both of the
sums $\partial(c)$ and $\partial(c')$, and with the same signs.  
The cells $c_{\tau(e_{i})}, c'_{\tau(e'_{i})}$  occur in the 
sums $\partial(c)$ and $\partial(c')$ (respectively), and with the same signs.  
\qed
\end{lemma}

\section{Background on the Cohomology Rings of Tree Braid
Groups}

\subsection{The cup product complex $\widehat{\ud{n}{T}}$}

In \cite{FS2}, we introduced a complex $\widehat{\ud{n}{T}}$ as an aid
to computing cup products in $H^{\ast} \left( \ud{n}{T} \right)$.  Here
we summarize the definition of this complex and collect a few useful facts.
We emphasize the standing assumption that $T$ is a tree, since certain
of the following propositions (\ref{ub}, for instance) aren't necessarily
true for general graphs.

If $c= \{ c_1 , c_2 , \ldots , c_n \}$ and 
$c' = \{ c_{1}' , c_{2}' , \ldots , c_{n}' \}$ are cells of $\ud{n}{T}$,
write $c \sim c'$ if:
\begin{enumerate}
\item $c$ and $c'$ have the same collections of edges, and
\item if $\{ e_1 , \ldots, e_k \}$ and $\{ e_{1}', \ldots e_{k}' \}$
are the edges in $c$ and $c'$ (respectively) then each component 
of 
$T - ( e_1 \cup e_2 \cup \ldots \cup e_k ) =        
T - ( e_{1}' \cup e_{2}' \cup \ldots \cup e_{k}' )$
contains the same number of vertices from $c$ as from $c'$.
\end{enumerate}
It is obvious that $\sim$ is an equivalence relation.  We let $[c]$ denote
the equivalence class of $c$.

There is a natural partial order $\leq$ defined on equivalence classes.
Write $[c] \leq [c']$ if there exist representatives $c \in [c]$, 
$c' \in [c']$ (without loss of generality) such that $c$ is a face of
$c'$ in the complex $\ud{n}{T}$.  In other words, $[c] \leq [c']$ if
there is a collection $\{ e_{i_1}' , e_{i_2}' , \ldots , e_{i_j}' \}
\subseteq c' $ such that $c$ may be obtained from $c'$ by replacing the 
edges $e_{i_1}', \ldots, e_{i_j}'$ in $c'$ with some combination of their initial and terminal
vertices.  

\begin{proposition} \label{ub} \cite{FS2} 
The relation $\leq$ is a well-defined partial order
with the property that any collection $\{ [c_1] , \ldots , [c_m] \}$
having an upper bound also has a least upper bound. 
\qed
\end{proposition}


\begin{proposition} \label{happy} \cite{FS2}      
If $\tilde{c}$ is a $j$-cell in $\ud{n}{T}$, then there is a unique
collection $\{ [c_1], \ldots, [c_j] \}$ of equivalence classes of
$1$-cells such that $[\tilde{c}]$ is the least upper bound of
$\{ [c_1], \ldots, [c_j] \}$.  

Indeed, $\{ [c_1], \ldots, [c_j] \} = \{ [c'] \mid [c'] \leq [\tilde{c}];
\, \, \, \mathrm{dim}(c')=1 \}$.
\qed
\end{proposition}

These  propositions 
establish a natural one-to-one correspondence between
equivalence classes of $j$-cells and $j$-element collections of equivalence
classes of $1$-cells having an upper bound.  

For
each equivalence class $[c]$ of 1-cells, introduce a copy of $S^{1}$,
denoted $S^{1}_{[c]}$, where $S^{1}_{[c]}$ is given the CW structure having
a single open $1$-cell $e^{1}_{[c]}$
 and a single $0$-cell.  We choose an orientation 
for each $1$-cell $e^{1}_{[c]}$
(or, in other words, an 
equivalence class of characteristic map 
$h_{[c]}: [0,1] \rightarrow \overline{e^{1}_{[c]}}$). 
The complex $\widehat{\ud{n}{T}}$
is the subcomplex of $\prod_{[c]} S^{1}_{[c]}$ 
obtained by throwing out
all open cells of the form 
$e^{1}_{[c_1]} \times \ldots \times e^{1}_{[c_k]}$ where the
collection $\{ [c_1] , \ldots , [c_k] \}$ has no upper bound.  If the
collection $\{ [c_1] , \ldots , [c_k] \}$ (where $c_{1}, \ldots, c_{k}$
are $1$-cells)
has a least upper bound $[c]$, 
then we label
the open cell $e^{1}_{[c_1]} \times \ldots \times e^{1}_{[c_k]}$   
by $[c]$.  Each cell of $\widehat{\udn{T}}$ is given the natural 
product orientation.  (Note that this orientation depends upon
an ordering of the factors.  We will order the factors so that
$[c_i ]$ occurs before $[c_j ]$ if $\iota(e_{i}) < \iota (e_{j})$, where
$e_{i}$ is the unique edge occurring in $c_i$, and $e_{j}$ is the unique edge
occurring in $c_j$.  If $\iota( e_{i}) = \iota( e_{j} )$, then the ordering may be chosen
arbitrarily.) 
Propositions \ref{ub} and \ref{happy} show that every equivalence
class $[c]$ of cells occurs as the label of some cell in 
$\widehat{\ud{n}{T}}$, and that each such label occurs uniquely. 

The complex $\widehat{\udn{T}}$ can also be described as  
a quotient of $\udn{T}$.  
We define a quotient map $q: \udn{T} \rightarrow \widehat{\udn{T}}$ as 
follows.     
The 
map $q$ sends the $j$-cell     
$c = \{ e_1 , e_2, \ldots , e_j , v_{j+1}, \ldots,  v_{n} \}$ in $\udn{T}$, where 
$\iota(e_1) < \iota(e_2) < \ldots < \iota(e_j)$,
to the cell
labelled $[c]$ in $\udn{T}$.  Let $f: [0,1]^{j} \rightarrow c$ be the standard characteristic
map for $c$ (as described at the end of Section \ref{gb}). For each cell $[c]$ in $\widehat{\udn{T}}$,
we choose an orientation $h: [0,1]^{j} \rightarrow [c]$
satisfying $h( t_1 , \ldots, t_j ) = (h_{[c_1]}(t_1) , \ldots, h_{[c_j]}(t_j))$, where $\{ [c_1], \ldots, [c_j] \}$
is the unique collection of equivalence classes of $1$-cells having $[c]$ as its least upper bound,
and $[c_i]$ is the unique equivalence class of this collection satisfying $e_i \in c_i$.  We define 
the map $q: \udn{T} \rightarrow \widehat{\udn{T}}$ cell by cell to be the unique map satisfying
$q \circ f = h$.  
The argument from 4.3 in \cite{FS2} shows that this assignment does indeed induce a quotient map
$q: \udn{T} \rightarrow \widehat{\udn{T}}$, and that each open cell $c$ in the domain is mapped homeomorphically
to the open cell $[c]$.  Our current choice of orientation also shows that positive cells are mapped to positive cells.   

We would now like to describe the cohomology rings $H^{\ast}( \widehat{\udn{T}}; \mathbb{Z})$.
For this we need a definition.  If $K$ is a finite simplicial complex, then the \emph{exterior face ring
of $K$} over $\mathbb{Z}$, denoted $\Lambda [ K ]$, is defined by 
$$ \Lambda [K] \cong \Lambda [ v_{1}, \ldots, v_{m} ] / I ,$$
where $\Lambda [ v_{1}, v_{2}, \ldots, v_{m} ]$ denotes the ordinary integral
exterior ring \cite{Hatcher} generated by the collection of vertices of $K$, and $I$ is the ideal
generated by the monomials $v_{i_{1}} v_{i_{2}} \ldots v_{i_{j}}$ such that
$\{ v_{i_{1}}, \ldots, v_{i_{j}} \}$ is not a simplex of $K$. 

It is now straightforward to describe the integral cohomology ring of
$\widehat{\ud{n}{T}}$.  For any equivalence class $[c]$ of
$j$-cells, let $\hat{\phi}_{[c]}$ denote the $j$-cocycle
satisfying $\hat{\phi}_{[c]}([c]) = 1$ and $\hat{\phi}_{[c]}([c']) =0$
for all $[c'] \neq [c]$.  We note that there is no real distinction
between cocycles and cohomology classes (since all boundary maps in the
cellular chain complex for $\widehat{\ud{n}{T}}$ are $0$), so we also
let $\hat{\phi}_{[c]}$ denote a cohomology class whenever 
it is convenient to do so. 

\begin{proposition} \label{mult} 
\cite{FS2}     
The integral cohomology ring $H^{\ast} ( \widehat{\ud{n}{T}} )$
is isomorphic to 
$ \Lambda \left[ K \right],$
where $K$ is a simplicial complex on the vertex set $V$ consisting of all
equivalence classes of $1$-cells in $\ud{n}{T}$. The simplices of $K$ consist
of all sets  
$\left\{ [c_{i_1}], \ldots, [c_{i_m}] \right\} \subseteq V$ 
having an upper bound.   

The isomorphism $\Phi: H^{\ast} ( \widehat{\ud{n}{T}})
\rightarrow \Lambda \left[ K \right]$ sends a $j$-dimensional 
cohomology class
$\hat{\phi}_{[c]}$ to $[c_1][c_2] \ldots [c_j]$, where 
$\{ [c_1], \ldots, [c_j] \}$ is the unique collection having $[c]$ as
its least upper bound. 

Here the ordering of the factors in the product $[c_1] [c_2] \ldots [c_j]$ is such that
$\iota(e_1) < \iota(e_2) < \ldots < \iota(e_j)$, where $e_i$ is the unique edge in $c_i$, 
for $i = 1, \ldots, j$.  
\qed
\end{proposition}

\subsection{Standard cocycles on $\ud{n}{T}$}

We have the following description of the map 
$q^{\ast}: H^{\ast}( \widehat{\udn{T}}; \mathbb{Z})
\rightarrow H^{\ast}(\udn{T}; \mathbb{Z})$ on cohomology.

\begin{proposition} \label{quot} 
The 
map $q^{\ast}$ 
sends $\hat{\phi}_{[c]}$ to $\phi_{[c]} \in C^{\ast}( \udn{T}; \mathbb{Z})$ on the
level of cellular cocycles, where

\begin{center}
\begin{tabular} {ll}
$\phi_{[c]}(\tilde{c})=1$ & if $\tilde{c} \sim c$ \\
$\phi_{[c]}(\tilde{c})=0$ & otherwise.
\end{tabular}
\end{center}
\end{proposition}
\begin{proof}
This follows from the description of the map $q: \udn{T} \rightarrow \widehat{\udn{T}}$, 
and from our choice of orientation as in the previous subsection.
\end{proof}
         
We call the cocycles $\phi_{[c]}$ \emph{standard}.  A standard 
cocycle $\phi_{[c]}$ is called \emph{critical} if the equivalence class
$[c]$ contains at least one critical cell. 

\begin{lemma} \cite{FS2}     
If $[c]$ contains at least one critical cell, then it contains
exactly one, and every other cell in $[c]$ is redundant.
\qed
\end{lemma}   

\begin{corollary} \label{surjective}
The collection of critical cocycles represent a basis for
$H^{\ast}( \ud{n}{T})$.  In particular, $q^{\ast}$ is surjective.  
\end{corollary}

\begin{proof}
Let $\phi_{[c]}$ be a critical $i$-cocycle, and suppose without loss of 
generality that $c$ is a critical $i$-cell.  
We evaluate $\phi_{[c]}$ on 
cycle representatives for a basis
of $H_{i}(\ud{n}{T})$.  Recall that one such collection of cycles 
$\{ \hat{c}_1, \ldots , \hat{c}_{j} \}$ corresponds
to the set $\{ c_1 , c_2 , c_3 , \ldots, c_j \}$ of all critical $i$-cells, 
where a given cycle $\hat{c}_{k}$ has the form $c_{k} +$ 
(collapsible cells).

Consider $\phi_{[c]}( \hat{c}_{k})$.  The support of $\phi_{[c]}$ contains
no collapsible cells, so $\phi_{[c]}( \hat{c}_{k} ) = \phi_{[c]} ( c_k )$.
Since the support of $\phi_{[c]}$ contains exactly one critical cell
(namely $c$), it follows that $\phi_{[c]} ( \hat{c}_{k} ) = 1$ if
$c = c_k$, and $\phi_{[c]} ( \hat{c}_{k}) = 0$ otherwise.

It follows from this that the cohomology classes
$\langle \phi_{[c_1]} \rangle , \ldots, \langle \phi_{[c_j]} \rangle$
map to the dual basis $\hat{c}_{1}^{\ast}, \ldots, \hat{c}_{j}^{\ast}$ in
$Hom( H_{i}( \ud{n}{T} ), \mathbb{Z})$ via the universal coefficient
isomorphism $H^{i}( \ud{n}{T}) \rightarrow Hom( H_{i}(\ud{n}{T}), 
\mathbb{Z})$.
\end{proof}   

\section{A Presentation for $H^{\ast}\left( B_n T \right)$ as a Ring}

\subsection{Pictures of standard cocycles and other necessary definitions}

Let $\phi_{[c]}$ be a standard cocycle on $\ud{n}{T}$.  If 
$c= \{ v_1 , \ldots, v_{j}, e_{j+1}, \ldots, e_{n} \}$, then by definition
$\phi_{[c]}$ is supported on all cells $\hat{c} = \{ \hat{v}_{1},
\ldots, \hat{v}_{j}, \hat{e}_{j+1}, \ldots, \hat{e}_{n} \}$ such that:
i) $\{ e_{j+1} , \ldots, e_{n} \} = 
\{ \hat{e}_{j+1}, \ldots, \hat{e}_{n} \}$ and 
ii) $| C \cap \{ v_1 , \ldots, v_{j} \} | = | C \cap \{ \hat{v}_{1}, 
\ldots, \hat{v}_{j} \} |$ for every connected component $C$ of 
$T - ( e_{j+1} \cup \ldots \cup e_{n} )$ ( $= T - ( \hat{e}_{j+1} \cup
\ldots \cup \hat{e}_{n} )$).  Thus, a standard $i$-cocycle is completely
determined by: i) a collection of connected components 
$Comp ( T - (  e_1 \cup \ldots \cup e_i ) )$ for a given collection
$E = \{ e_1 , \ldots, e_i \}$ of disjoint closed edges,  
and ii)
an assignment $f: Comp( T - ( e_1 \cup \ldots \cup e_{i} )) \rightarrow
\mathbb{Z}^{+} \cup \{ 0 \}$ such that $\sum_{c} f(c) = n-i$.  We call
a pair $( Comp(T- (e_1 \cup \ldots \cup e_{i}) ), f)$ satisfying
i) and ii) a \emph{cloud picture} of a standard cocycle, since
we will often indicate such a pair using pictures as in the following
example.  (Note:  in practice, 
we will usually think of the elements of $Comp ( T - ( e_1 \cup \ldots \cup e_i ))$
as the $0$-skeletons of the connected components described above in i).)  

\begin{example}
Let $T$ be as in Figure \ref{number}, and consider pairs $( Comp( T- e_{19}), f_1 )$,
$( Comp(T- e_{16} ), f_2)$ (where $e_i$ is the unique edge of $T$ satisfying
$\tau(e_i) = v_{i}$, for $i = 16,19$), and 
$f_1 : \{ C_0 , C_1 , C_2 \} \rightarrow \mathbb{Z}^{+} \cup \{ 0 \}$,
$f_2 : \{ C'_0 , C'_1 , C'_2 \} \rightarrow \mathbb{Z}^{+} \cup \{ 0 \}$
satisfying:
\begin{eqnarray*}
f_{1}(C_0 ) = 1 & ; & f_{2}(C'_0) = 2 \\
f_{1}(C_1 ) = 2 & ; & f_{2}(C'_1) = 1 \\
f_{1}(C_2 ) = 0 & ; & f_{2}(C'_2) = 0.   
\end{eqnarray*}
(Here $C_i$ is the connected component of $T - e_{19}$ lying in the direction
$i$ from $\iota( e_{19} ) = v_{9}$; similarly $C'_{i}$ lies in direction
$i$ from $\iota( e_{16} )$.)

We represent each of these pairs as pictures (see Figure \ref{cup1}).
Each picture can be interpreted in two ways.  First, we can see each picture
as a description of a given equivalence class of, in this case, $1$-cells.
The picture on the left, for instance, represents the equivalence class
of the $1$-cell $c = \{ e_{19} , \ast , v_{10}, v_{11} \}$, where $\ast$ is the basepoint, and 
$v_{10}$ and $v_{11}$ are the vertices numbered $10$ and $11$ in Figure \ref{number}.  
The elements of the indicated equivalence class are $1$-cells $c'$ such that: i) the unique (in this case)
$1$-cell of $c'$ is $e_{19} $, and ii) exactly two vertices of $c'$ lie inside the cloud at the bottom left
(which consists of the set $\{ v_{10}, \ldots, v_{18} \}$), iii) exactly one vertex of $c'$ lies inside 
the cloud at the bottom right (which consists of the vertices $\{ \ast, v_{1}, \ldots, v_8 \}$).  We leave
it as an easy exercise to describe the equivalence class symbolized by the picture on the right in Figure 
\ref{cup1}.  Second, we can see each picture as a representative of
the unique standard cocycle supported on the given equivalence class. 

\begin{figure} [!h]
\begin{center}
\includegraphics{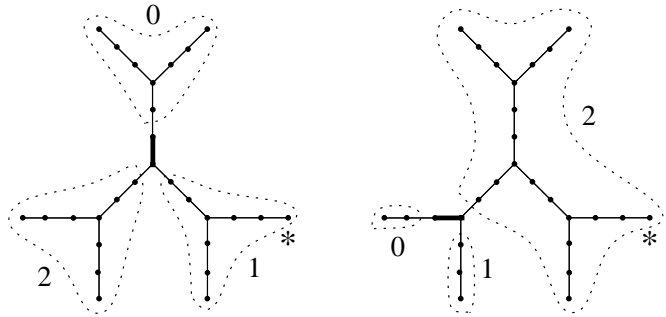}
\caption{These cloud pictures are determined by the
pairs $( Comp( T- e_{19}), f_1 )$ and $( Comp( T- e_{16}), f_2 )$, respectively.} 
\label{cup1}
\end{center}
\end{figure}

Now Propositions \ref{mult} and \ref{quot} imply that the product of the above
two cocycles (cohomology classes), in the given order, is as in Figure \ref{cup2}.    
The main point to check is that the cloud picture in Figure \ref{cup2} is the least upper bound
of the two pictures in Figure \ref{cup1}.  This can be done very easily.  Simply choose a representative $c$
of the equivalence class $[c]$ determined by the picture in Figure \ref{cup2}, and find two faces $c_1$, $c_2$ of
$c$ such that $[c_1]$ and $[c_2]$ have the pictures in Figure \ref{cup1}.  One possibility is
to let $c = \{ e_{19} , e_{16}, \ast, v_{13} \}$, $c_1 = \{  e_{19} , v_{12}, \ast , v_{13} \}$, 
and $c_2 = \{ v_{9}, e_{16}, \ast, v_{13} \}$.       

\begin{figure} [!h]
\begin{center}
\includegraphics{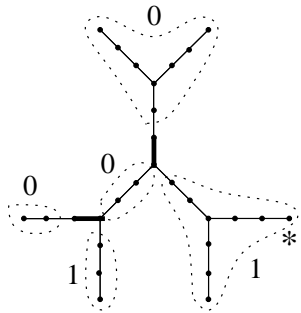}
\caption{This is the product of the two cocycles from Figure \ref{cup1}.}
\label{cup2}
\end{center}
\end{figure}

We introduce some definitions which will be useful when we work with standard cocycles.  First, we refer to 
the elements $C$ of $Comp (T - (e_1 \cup \ldots \cup e_j ))$ as \emph{clouds} of the given standard cocycle. 
If $\mathcal{C}$
is a cloud picture and $e$ is an edge of $\mathcal{C}$, we say that the cloud picture $\mathcal{C}'$ is
obtained from $\mathcal{C}$ by \emph{breaking the edge $e$ and combining clouds} if $\mathcal{C}'$ is the
result of:  i) removing the edge $e$ from the defining list for $\mathcal{C}$, and ii) replacing the clouds
incident with $e$ by a single cloud (which also contains the edge $e$).  The resulting cloud's number is $1$ 
greater than the sum of the numberings of the old clouds.  For instance, the cloud picture on the left in Figure
\ref{cup1} is obtained from the cloud picture in Figure \ref{cup2} by breaking the edge $e_{16}$ and combining
clouds.  

\qed
\end{example}

We now assign an ordered pair of integers to each standard cocycle.  We need
several definitions.  
If $c= \{ e_{1}, \ldots, e_{i}, v_{i+1}, \ldots, v_{n} \}$ is an $i$-cell
of $\ud{n}{T}$, then we call the vertices $\iota(e_1 ), \ldots, \iota(e_{i})$
\emph{special} vertices of $c$. The \emph{$c$-depth} $dep_{c} (v_k)$ of
a vertex $v_k \in c$ is the number of special vertices lying on the geodesic
segment $[ \ast, v_{k} ]$ in $T$.  The \emph{total vertex depth} $dep(c)$ of
$c$ is $\sum_{v \in c} dep_{c}(v)$.

If $c$ is an arbitrary cell
of $\ud{n}{T}$, let $r(c)$ be a choice of cell in $[c]$ 
with the property that all of its vertices are blocked.  
An edge
$e \in c$ is called \emph{bad} if it is order-respecting in $r(c)$.  Note
that ``badness" of an edge is well-defined on equivalence classes $[c]$,
in the sense that if $c \sim c'$ and $e \in c,c'$ then $e$ is bad in $c$
if and only if $e$ is bad in $c'$.  (On the other hand, easy examples show
that there exist $c$, $c'$, $c \sim c'$ for which it is possible to choose
$r(c)$ and $r(c')$ such that $r(c) \neq r(c')$.)  The
number of bad edges in $c$ is denoted $bad(c)$.

We define the \emph{rank of $c$}, denoted $rank(c)$, to be the pair
$( - dep(c), bad(c) )$.  The set of ranks is lexicographically ordered,
i.e., $rank(c) < rank(c')$ if $dep(c) > dep(c')$ or $dep(c) = dep(c')$ and
$bad(c) < bad(c')$.  It is clear moreover that there can be no infinite 
descending chains.  Note that the functions $dep$ and $bad$ are
well-defined on equivalence classes, so we can define the rank $rank([c])$
of an equivalence class to be equal to the rank of $c$.  Similarly, 
the \emph{rank} of a standard cocycle is the rank of any cell in its support.

\subsection{A complete list of relators}

Let $\mathcal{C} = ( Comp(T - (e_1 \cup \ldots \cup e_i )), f)$ be a
cloud picture of a standard $i$-cocycle; let $E = \{ e_1 , \ldots, e_{i} \}$,
and $Comp(T - (e_{1} \cup \ldots \cup e_{i} )) = \{ C_1 , \ldots , C_k \}$.
Let $v'$ be a special vertex of one (equivalently, any) of the cells $c$
in the support of the cocycle determined by $\mathcal{C}$.  Let $e'$ be 
the unique edge of $c$ having $v'$ as its initial endpoint.  We define
an $(i-1)$-dimensional cochain 
$R_{\mathcal{C},v'} =  \left( \mathcal{P}_{\mathcal{C}, v'}, 
f_{\mathcal{C}, v'} \right)$ in terms of a modified ``cloud notation". 
The set $\mathcal{P}_{\mathcal{C}, v'}$ partitions the vertices
of $T - \left( \cup_{e \in E - \{ e' \}} e \right)$; each member of 
$\mathcal{P}_{\mathcal{C}, v'}$ is precisely the same as a member of
$Comp( T- ( e_1 \cup \ldots \cup e_{i} ))$, with two exceptions:  i)
the component $C_{\tau} \in Comp( T - (e_1 \cup \ldots \cup e_{i}))$ which
is adjacent to $\tau(e')$ is replaced by the set 
$C'_{\tau} = C_{\tau} \cup \{ \tau(e') \}$ in $\mathcal{P}_{\mathcal{C},v'}$,
and ii) the component $C_{\iota} \in Comp( T - (e_1 \cup \ldots \cup e_{i}))$
which is adjacent to $\iota(e') = v'$ and lies in the $0$th direction from
$v'$ is replaced by the set $C'_{\iota} = C_{\iota} \cup \{ \iota(e') \}$.
The function $f_{\mathcal{C},v'}$: i) agrees with $f$ on 
$\mathcal{P}_{\mathcal{C}, v'} \cap 
Comp( T - ( e_{1} \cup \ldots \cup e_{i} ))$, 
ii) $f_{\mathcal{C},v'}(C'_{\iota}) = f( C_{\iota})$, and 
iii) $f_{\mathcal{C},v'}(C'_{\tau})   
= f(C_{\tau}) +1$.  

This new ``cloud notation" describes the support of a 
cochain $R_{\mathcal{C},v'}$.  It is perhaps best to give the definition
of this cochain with help from an example (Example \ref{basicex} below).

Note that our definition of $R_{\mathcal{C}, v'}$ implicitly assumes
that the degree of the vertex $\tau(e')$ is less than or equal to $2$.
This is by far the most important case.  We will treat the case
in which the degree of $\tau(e')$ is at least $3$ after Theorem \ref{biggie}.  

We now describe such a cochain $R_{\mathcal{C},v'}$ and its coboundary
$\delta( R_{\mathcal{C},v'})$ in an explicit case.

\begin{example} \label{basicex}
Consider the standard cocycle on the left in Figure \ref{mainex}.
\begin{figure}[!h]
\begin{center}
\includegraphics{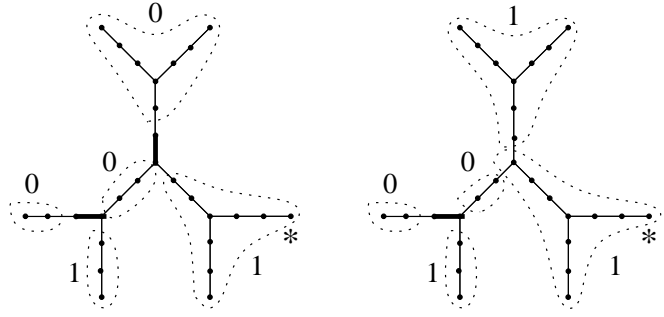}
\caption{Here is a standard cocycle $\phi$ and an associated cochain $R_{\mathcal{C}, v'}$
whose coboundary can be used to express $\phi$ in terms of standard cocycles of smaller rank.}
\label{mainex}
\end{center}
\end{figure}
If we check the above definitions, we see that the $1$-dimensional cochain
$R_{\mathcal{C}, v_9}$ is pictured on the right in ``cloud notation".  
Thus a positively oriented cell $c$ of
$\ud{n}{T}$ lies in the support of $R_{\mathcal{C}, v_9}$ ( and $R_{\mathcal{C}, v_9} (c) = 1$)
if and only if:
i) the only edge of $c$ is $e_{16}$; ii) $c$ contains exactly one of the
cells from the set $\{ v_{13}, v_{14}, v_{15} \}$; iii) $c$ contains
exactly one of the cells from $\{ \ast, v_1 , \ldots, v_{9} \}$, and iv)
$c$ contains exactly one of the cells from $\{ v_{19}, \ldots, v_{27} \}$.

Let us consider the support of the ($2$-dimensional) coboundary
$\delta(R_{\mathcal{C},v_9})$.  Suppose that $c' \in supp \, \delta ( R_{\mathcal{C},v_9})$.
Let $e'_{1}, e_{2}' \in c'$ be the two 
$1$-cells in $c' \in supp \, \delta( R_{\mathcal{C},v_9})$.  
The coboundary $\delta( R_{\mathcal{C}, v_9})$ evaluated at $c'$ is the
following sum for a suitable choice of sign:
$$ \delta( R_{\mathcal{C}, v_9})(c') = \pm \left(
R_{\mathcal{C},v_9} ( c'_{\tau(e'_{1})}) -
R_{\mathcal{C},v_9} ( c'_{\iota(e'_{1})}) \right) \mp \left( 
R_{\mathcal{C},v_9} ( c'_{\tau(e'_{2})}) -
R_{\mathcal{C},v_9} ( c'_{\iota(e'_{2})}) \right),$$
where, for instance, $c'_{\tau(e'_{1})}$ is the face of $c'$ obtained
by replacing $e'_{1} \in c'$ with $\tau(e'_{1})$.  The cells 
$c'_{\iota(e'_{1})}$,
$c'_{\tau(e'_{2})}$, and 
$c'_{\iota(e'_{2})}$ have analogous definitions.

It follows easily from this that one of $e'_{1}$, $e'_{2}$ is $e_{16}$
(since $c' \in supp \, \delta(R_{\mathcal{C}, v_9})$); we assume
$e'_{1} = e_{16}$, without loss of generality.  The cells 
$c'_{\tau(e'_{1})},   
c'_{\iota(e'_{1})}$ 
are thus annihilated by $R_{\mathcal{C},v_9}$.
We have
$$ \delta( R_{\mathcal{C}, v_9})(c') = 
\mp \left( 
R_{\mathcal{C},v_9} ( c'_{\tau(e'_{2})}) -
R_{\mathcal{C},v_9} ( c'_{\iota(e'_{2})}) \right).$$
It follows directly from this that precisely one of
$c'_{\tau(e'_{2})},   
c'_{\iota(e'_{2})}$ lies in the support of $R_{\mathcal{C},v_9}$.  It
then follows that $e'_{2}$ must connect a vertex in one ``cloud" to a
vertex in another.  Thus $e'_{2}$ is either $e_{10}$ or $e_{19}$.
(Note also that if $e'_{2}$ connects two clouds as above, then it follows
easily that at most one face of $c'$ lies in the support of $R_{\mathcal{C}, v_9}$.)

The support of $\delta(R_{\mathcal{C}, v_9})$ can now be 
determined directly.  If $c' \in supp \, \delta( R_{\mathcal{C}, v_9})$,
then 
$c' = \{ e_{16}, e_{10}, v' , v'' \}$ or
$c' = \{ e_{16}, e_{19}, v', v'' \}$.  
The above reasoning also shows
that the entire contribution to the value of $\delta(R_{\mathcal{C}, v_9 })$
comes from a single face $c''$ of $c'$ lying in the support of
$R_{\mathcal{C}, v_9}$.  There are four cases: 
\begin{enumerate}
\item 
$c' = \{ e_{16}, e_{10}, v' , v'' \}$ and   
$c'' = \{ e_{16}, \iota(e_{10}), v' , v'' \}$; 
\item 
$c' = \{ e_{16}, e_{10}, v' , v'' \}$ and   
$c'' = \{ e_{16}, \tau(e_{10}), v' , v'' \}$; 
\item 
$c' = \{ e_{16}, e_{19}, v', v'' \}$ and    
$c'' = \{ e_{16}, \iota(e_{19}), v' , v'' \}$; 
\item 
$c' = \{ e_{16}, e_{19}, v', v'' \}$ and   
$c'' = \{ e_{16}, \tau(e_{19}), v' , v'' \}$.
\end{enumerate}

\begin{figure} [!h]
\begin{center}
\includegraphics{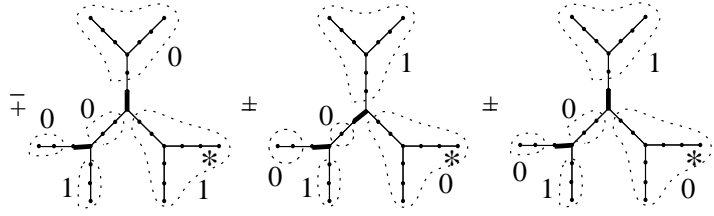}
\caption{This is the coboundary $\delta( R_{\mathcal{C}, v_9 })$.}
\label{mainex2}
\end{center}
\end{figure}

Case (1) can be realized    
if $v' \in \{ v_{19}, \ldots, v_{27} \}$ and $v'' \in 
\{ v_{13}, v_{14}, v_{15} \}$.  Case (2) is clearly impossible, since
$\tau( e_{10})$ is in a cloud which is labelled $0$, and thus
$c''$ fails to be in $supp( R_{\mathcal{C}, v_9 })$.  Case (3) can be
realized if 
$v' \in \{ v_{20}, \ldots , v_{27} \}$ and $v'' \in 
\{ v_{13}, v_{14}, v_{15}  \}$.  Case (4) can be realized 
if $v' \in \{ \ast,
v_{1}, \ldots, v_{8} \}$ and $v'' \in \{ v_{13}, v_{14}, v_{15} \}$.

Reviewing the considerations in (1)-(4) above, we see that
$supp \, \delta( R_{\mathcal{C}, v_9})$ is the union of the (disjoint)
supports of three standard cocycles; their cloud pictures are 
Figure \ref{mainex2}, written (from left to right) in the order 
(4), (1), (3).
The sum in Figure \ref{mainex2} 
is equal to $0$ on cohomology since it is a coboundary. 
Note that the standard cocycle on the left is $\phi$. 

Now $rank(\phi) = (-2,1)$;
the ranks of the standard cocycles on the right are both $(-3,1)$.  We therefore
have succeeded in representing the cohomology class of $\phi$  
using a sum of cocycles of rank smaller than $rank(\phi)$. \qed
\end{example} 

We now consider the features from the previous example which hold more
generally.  Let $\mathcal{C}$ be the cloud picture of a standard cocycle
$\phi$ (of any dimension), let $v$ be a special vertex of $\phi$, and let
$e$ be the (unique) edge of $\phi$ such that $\iota(e) = v$.  A cell
$c'$ in the support of $\delta \left( R_{\mathcal{C}, v} \right)$ has
the following properties:
\begin{enumerate}
\item The edge set $E( c')$ contains the edge set 
$E \left( R_{\mathcal{C}, v} \right)$.  The one remaining edge
$e' \in E(c')$ must connect distinct clouds of $R_{\mathcal{C},v}$.  
Thus, there are (a priori, at least) $d(v) -1$ possibilities for $e'$.
\item Exactly one face of $c'$ lies in the support of
$R_{\mathcal{C},v}$, and this face can be obtained by replacing the
edge $e'$ with either its initial or its terminal endpoint.  We therefore
have the formula
$$ \delta \left( R_{\mathcal{C},v} \right) (c') = 
R_{\mathcal{C},v} \left( c'_{\tau(e')} \right) - 
R_{\mathcal{C},v} \left( c'_{\iota(e')} \right),$$
where, for example, $c'_{\tau(e')}$ is the face of $c'$ obtained by replacing
the edge $e'$ with $\tau(e')$ in the definition of $c'$.  This equality
holds only up to a choice of multiplication of the right side of the
equality by $\pm 1$, but this choice doesn't matter, since we are only
interested in equalities on the level of cohomology, and the left side of the equation
is, of course, $0$ on cohomology.  We ignore this choice
of sign from now on.
\item We can divide the cells $c'$ into two groups, according to whether
$c'_{\tau(e')}$ or $c'_{\iota(e')}$ lies in the support of 
$R_{\mathcal{C},v}$.  By Lemma \ref{miracle}, we can write
$$ \delta \left( R_{\mathcal{C},v} \right) = \Sigma_{\mathcal{C}, v, \tau}
- \Sigma_{\mathcal{C}, v, \iota}.$$
Here 
$$ \Sigma_{\mathcal{C}, v, \tau} = \chi \left( \left\{ c' \mid
c'_{\tau(e')} \in supp \, R_{\mathcal{C},v} \right\} \right)$$
and $\chi$ denotes a characteristic function.  The function
$\Sigma_{\mathcal{C}, v, \iota}$ is defined similarly.
\item We now examine the functions $\Sigma_{\mathcal{C}, v, \tau}$
and $\Sigma_{\mathcal{C}, v, \iota}$.  We can partition the support
of $\Sigma_{\mathcal{C}, v, \tau}$ into $d(v)-1$ pieces, according to
the location of the edge $e'$, which can point in any of $d(v)-1$ directions.
(In practice, some elements of the indicated partition will be empty, 
as in Case (2) from Example \ref{basicex}, but
this causes no problems.  It will simply mean that some of the standard
cocycles in the sums to be defined below are $0$.)
For $i=1, \ldots, d(v)-1$, we let $e'_{i}$ denote the edge
with initial vertex $v$ which lies in the direction $i$ from $v$.  For
$i = 0, 1, \ldots, d(v)-1$, we let $C_{i}$ denote the element of the
partition $\mathcal{P}_{\mathcal{C}, v}$ which lies in the $i$th direction
from $v$.  (Recall that $R_{\mathcal{C},v} = ( \mathcal{P}_{\mathcal{C},v},
f_{\mathcal{C},v})$.)

A routine check shows that
$$ \Sigma_{\mathcal{C}, v, \tau} = \sum_{i=1}^{d(v)-1} 
\Psi_{\mathcal{C}, v, \tau, i}$$
where $\Psi_{\mathcal{C}, v, \tau, i}$ is the standard cocycle supported
on all cells $c$ such that: i) $E(c) = E( R_{\mathcal{C},v}) \cup
\{ e'_{i} \}$, and ii) exactly $f_{\mathcal{C}, v} ( c_j )$ vertices of
$c$ lie in $C_j$ for $j \neq i$; exactly 
$f_{\mathcal{C}, v} ( C_i ) - 1$ vertices of $c$ lie in $C_i$.

Similarly,
$$ \Sigma_{\mathcal{C}, v, \iota} = \sum_{i=1}^{d(v)-1}
\Psi_{\mathcal{C}, v, \iota, i}$$
where $\Psi_{\mathcal{C}, v, \iota, i}$ is the standard cocycle supported
on all cells $c$ such that:  i) 
$E(c) = E( R_{\mathcal{C}, v}) \cup \{ e'_{i} \}$, and
ii) exactly $f_{\mathcal{C}, v} (C_{i})$ vertices of $c$ lie in $C_i$ for $i>0$, 
and $f_{\mathcal{C}, v} (C_{0}) - 1$ vertices of $c$ lie in $C_{0}$. 
\end{enumerate}

Since the descriptions of $\Sigma_{\mathcal{C}, v, \tau}$ and
$\Sigma_{\mathcal{C}, v, \iota}$ from (4) are somewhat unpleasant
notationally, we give a computation of $\delta( R_{\mathcal{C}, v})
= \Sigma_{\mathcal{C}, v, \tau} - \Sigma_{\mathcal{C}, v, \iota}$ in the
example below.

\begin{example}
Let $R_{\mathcal{C},v}$ denote the cochain determined by the cloud
representative in Figure \ref{ex} a).  We have abandoned the cloud notation, 
\begin{figure}[!h]
\begin{center}
\includegraphics{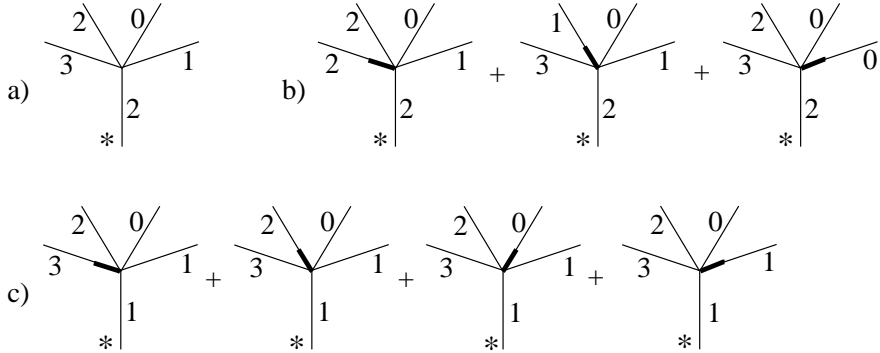}
\caption{a) shows a picture representative of the cochain $R_{\mathcal{C},v}$; b) is the
sum $\Sigma_{\mathcal{C}, v, \tau}$, and c) is $\Sigma_{\mathcal{C}, v, \iota}$.}
\label{ex}
\end{center}
\end{figure}
omitted the numbering, and even avoided specifying the subdivision of the tree $T$.
(A sufficient subdivision of $T$ would have $7$ edges along each of the radial arms
of the tree.)  In a), the integer labels on the arms of the tree represent the
number of vertices which lie inside the cloud containing the given arm.  For this
purpose, we regard the central vertex as belonging to the cloud which also contains the
basepoint.  The notation in b) and c) should be unambiguous.

The two constituent sums of the coboundary $\delta( R_{\mathcal{C},v} )$ are
depicted in Figure \ref{ex}, parts b) and c).
The terms in b) are $\Sigma_{\mathcal{C}, v, \tau}$;
the terms in c) are 
$\Sigma_{\mathcal{C}, v, \iota}$. \qed
\end{example}

Note that the ranks of the cocycles in Figure \ref{ex}, parts b) and c), are
$(-5,1), (-5,0), (-5,0), (-6,1), (-6,0), (-6,0),$ and $(-6,0)$, respectively
(read from left to right).  Thus there is a unique standard cocycle in the
sum $\Sigma_{\mathcal{C}, v, \tau} - \Sigma_{\mathcal{C}, v, \iota}$ having
maximal rank (namely the first cocycle, in this case).  This is in fact a
general feature of all coboundaries $\delta( R_{\mathcal{C}, v})$ with
$\mathcal{C}$ and $v$ arbitrary.

\begin{proposition} \label{shfp}
Let $\Psi$ be a standard cocycle with cloud diagram $\mathcal{C}$.  Let
$e$ be a bad edge of $\Psi$ and let $v$ be the initial vertex of
$e$.

The coboundary $\delta ( R_{\mathcal{C}, v})$ can be expressed as a sum
$\sum \Psi_{i}$ of standard cocycles, where
\begin{enumerate}
\item the cocycle $\Psi$ occurs exactly once in $\sum \Psi_{i}$, and
\item every other term in $\sum \Psi_{i}$ (if any) is a standard cocycle
of rank strictly less than $rank( \Psi )$.
\end{enumerate}
\end{proposition}
\begin{proof}
Suppose that $e$ lies in direction $k$ from $v$, where $k \in \{ 1, \ldots,
d(v)-1 \}$, or $k=1$ if $d(v) =1$.  Let $C_j$, for $j = 0, 1, \ldots, d(v)-1$,
be the cloud of $R_{\mathcal{C},v}$ lying in the $j$th direction from $v$.
It follows from the fact that $e$ is bad that 
$\phi_{\mathcal{C},v} ( C_j ) = 0$ for all $0 < j < k$; by the definition
of $R_{\mathcal{C},v}$, 
$f_{\mathcal{C},v} ( C_k ) = 1 + f_{\mathcal{C}}( C_k )$.  In particular,
$f_{\mathcal{C},v} (C_k ) > 0$.

Now we consider the sum $\Sigma_{\mathcal{C}, v, \tau}$:
\begin{eqnarray*}
\Sigma_{\mathcal{C}, v, \tau} 
& = & \sum_{i=1}^{d(v)-1} \Psi_{\mathcal{C},v, \tau, i} \\
& = & \Psi_{\mathcal{C}, v, \tau, k} + \sum_{i=k+1}^{d(v)-1}
\Psi_{\mathcal{C}, v, \tau, i}.
\end{eqnarray*}
The cocycle $\Psi_{\mathcal{C}, v, \tau, k}$ is equal to $\Psi$ by definition.
Each term of the remaining sum is either identically $0$ or the edge
$e'_{i}$ ($i>k$) is no longer bad.  Since the total vertex depth of each 
term in $\Sigma_{\mathcal{C}, v, \tau}$ is identical, and
$bad( \Psi_{\mathcal{C}, v, \tau, i}) < 
bad( \Psi_{\mathcal{C}, v, \tau, k})$
for $i>k$, it follows that $rank( \Psi_{\mathcal{C}, v, \tau, k} )
> rank( \Psi_{\mathcal{C}, v, \tau, i})$, for $i>k$.

In the sum
$$ \Sigma_{\mathcal{C}, v, \iota} = 
\sum_{i=1}^{d(v)-1} \Psi_{\mathcal{C}, v, \iota, i},$$
each cocycle $\Psi_{\mathcal{C}, v, \iota, i}$ satisfies 
$dep( \Psi_{\mathcal{C}, v, \iota, i}) = dep( \Psi ) + 1$, so
$rank( \Psi_{\mathcal{C}, v, \iota, i}) < rank( \Psi)$ for all $i$.
\end{proof}

Fix a cochain $R_{\mathcal{C},v}$, where $\mathcal{C}$ is a cloud picture and
$v$ is a special vertex of a bad edge of $\mathcal{C}$.  The discussion
before Proposition \ref{shfp} shows that $\delta ( R_{\mathcal{C},v} )$ 
can be written uniquely as a sum of standard cocycles.
(Indeed, the non-zero standard cocycles make up a linearly independent set, and even
have pairwise disjoint supports.)  We let $\widehat{ \delta ( R_{\mathcal{C},v})}$
denote the unique cellular cochain in $C^{\ast} ( \widehat{\udn{T}}; \mathbb{Z})$ which
maps onto $\delta( R_{\mathcal{C},v})$ under the function on cochains induced by the 
canonical quotient map
$q: \udn{T} \rightarrow \widehat{\udn{T}}$.  Since every cochain in $C^{\ast}( \widehat{\udn{T}}; \mathbb{Z})$
is a cocycle, we can safely identify $\widehat{ \delta( R_{\mathcal{C},v})}$ with the element of
cohomology it represents, and we shall do so freely and without further notice.

\begin{theorem} \label{biggie}
The kernel of the map 
$q^{\ast} : H^{\ast}( \widehat{\udn{T}}; \mathbb{Z} ) \rightarrow H^{\ast} ( \udn{T}; \mathbb{Z} )$
is generated by the cohomology classes $\widehat{ \delta( R_{\mathcal{C}, v})}$, where $\mathcal{C}$
is a cloud picture with at least one bad edge $e$ and $v$ is the initial vertex of $e$.

In particular,
$$ H^{\ast}( \udn{T} ; \mathbb{Z}) \cong \Lambda [K] / I, $$
where $I$ is the ideal generated by the $\widehat{ \delta( R_{\mathcal{C}, v})}$ (as above), and $K$
is the simplicial complex described in Proposition \ref{mult}.   
\end{theorem}

\begin{proof}
We work with standard cocycles in $C^{\ast}( \udn{T})$, which are in perfect one-to-one
correspondence with the canonical basis in $C^{\ast}( \widehat{\udn{T}})$ ($= H^{\ast}( \widehat{\udn{T}})$).

Define a strict partial order $<$ on standard cocycles as follows.  Let $\phi$ be a standard cocycle,
let $\mathcal{C}$ be its cloud picture, suppose that $\mathcal{C}$ has at least one bad edge $e$, and let $v= \iota(e)$.  
If $\phi_{1}$ is a standard cocycle we write $\phi_{1} < \phi$ if $supp \, \phi_{1} \subseteq supp \, \delta( R_{\mathcal{C},v})$
and $\phi_{1} \neq \phi$.   
(Thus, for instance, the left-most standard cocycle $\phi$ in Figure \ref{mainex2} is strictly greater than the other
cocycles in Figure \ref{mainex2}.) 
Let $<$ also denote the transitive closure of $<$.  That $<$ is 
a strict partial order follows easily from the fact that $rank( \phi_{1} ) < rank( \phi)$ by Proposition \ref{shfp}.

Consider a minimal standard cocycle $\phi$ under the partial order $<$.  By the definition of $<$, this implies either:
i) that $\phi$ has no bad edges, or ii) any sum $\delta( R_{\mathcal{C},v})$ (where $\mathcal{C}$ is the cloud picture
of $\phi$ and $v$ is the initial vertex of a bad edge) satisfies $\phi = \delta (R_{\mathcal{C},v})$.  In case i),
$\phi$ is a critical cocycle, and in case ii) $\phi$ is cohomologous to $0$.

Let $\sum_{i} \phi_{i}$ be a finite sum of standard cocycles.  We let the \emph{rank} of $\sum_{i} \phi_{i}$ be the maximum
of the ranks of the non-critical standard cocycles occurring in $\sum_{i} \phi_{i}$.  The rank of a sum of critical cocycles 
can simply be defined as $( - \infty , -\infty )$.  Suppose, for the sake of simplicity in notation, 
that there is a unique non-critical standard cocycle $\phi_{k}$ occurring in the sum such that $rank( \phi_{k}) > rank( \phi_{i})$
for any non-critical standard cocycle $\phi_{i}$ in $\sum_{i} \phi_{i}$ such that $i \neq k$.  The cloud picture $\mathcal{C}_{k}$ for
$\phi_{k}$ has at least one bad edge $e$, and we let $\iota(e)=v$.  We have the equality
$$ \sum_{i} \phi_{i} = \left( \sum_{i \neq k} \phi_{i} \right) + \left( \phi_{k} - \delta( R_{\mathcal{C}_{k}, v}) \right)$$
on the level of cohomology.  Note that, by Proposition \ref{shfp}, the cocycle on the right side is of smaller rank.

In the general case, there are several non-critical standard cocycles $\phi_{i_{1}}, \ldots, \phi_{i_{k}}$ (rather than a unique one) 
satisfying $rank(\phi_{i_{1}}) = rank ( \phi_{i_{2}} ) = \ldots = rank( \phi_{i_{k}} ) = rank ( \sum_{i} \phi_{i} )$.
In this case, the same procedure should be performed for each such $\phi_{i_{j}}$ simultaneously.

It follows that a given sum $\Sigma$ of standard cocycles, not all of which are critical, can be rewritten as $\Sigma'$, where
$rank( \Sigma' ) < rank( \Sigma )$.  The collection of all possible ranks is finite, so this rewriting procedure must terminate, 
and the result will be a (possibly empty) sum of critical cocycles.  During this process, we apply only relations (coboundaries) of
the form $\delta( R_{\mathcal{C},v})$ as in the statement of the theorem.
If the sum $\Sigma$ is cohomologous to $0$, then the rewriting procedure expresses $\Sigma$ as a sum of coboundaries $\delta (R_{\mathcal{C},v})$,
as required.
\end{proof}

Finally, we mention how to modify the above rewriting procedure to cover the case in which 
the degree of $\tau(e)$ is greater than $2$, for some defining edges $e$ of the standard cocycle $\phi$.
Fix such an edge $e$.  Consider the cochain $R'$, defined using  
``cloud notation", and satisfying:  i)  one of the clouds of $R'$ consists of the vertex $\tau(e)$ alone, 
and this cloud is assigned the number $1$; ii) each of the clouds incident with the edge $e$ is the same
as before, except for the cloud $C_0$ lying in direction $0$ from $\tau(e)$, which gains the vertex 
$\iota(e)$.  These clouds are all numbered as they were in $\phi$.

We can again express the coboundary $\delta( R' )$ as a sum $\Sigma_{\tau} - \Sigma_{\iota}$.  
We leave it as an exercise to show that $\Sigma_{\tau}$ consists simply of the standard cocycle $\phi$
itself, and that $\Sigma_{\iota}$ is a sum of standard cocycles in which the edge $e$ has been replaced by
various edges $e'$ (one for each term of the sum) satisfying $\iota( e' ) = \tau(e)$.  This in particular
implies that $\tau(e')$ has degree less than $3$ in every case.  (Here we assume that no two vertices of
degree greater than two are adjacent.  This can be guaranteed by subdividing $T$.)     

Now note that, after applying this procedure repeatedly, we can represent the cohomology class of $\phi$ by 
a sum of standard cocycles, none of which have edges $e$ such that $d( \tau(e) ) > 2$.  At this point, we
can apply the procedure of Theorem \ref{biggie} without change, and ultimately express the cohomology class of $\phi$
in terms of critical cocycles.  Note that the procedure of that Theorem will never change the special vertices
of a standard cocycle, so it is never necessary to apply the special procedure described here again.  

\subsection{A Refinement of Ghrist's Conjecture} \label{ghrist}

The main theorem and (especially) the examples to come in Section \ref{exples} suggest a
possible general description of the cohomology rings $H^{\ast}( \udn{T}; \mathbb{Z})$.

\begin{conjecture} \label{newconj}
Fix a finite tree $T$ and a positive integer $n$.  Let $K$ be the simplicial complex
associated to the complex $\widehat{\udn{T}}$, as in Proposition \ref{mult}.
$$ H^{\ast}( \udn{T}; \mathbb{Z} ) \cong \Lambda[K] / I',$$
where $I'$ is the ideal generated by all monomials $[c_{1}]\ldots[c_{k}]$, where the
least upper bound $[c]$ of $\{ [c_{1}], \ldots, [c_{k}] \}$ contains no critical cell,
where we regard $[c]$ as an equivalence class of cells in $\udn{T}$.

In particular, each $H^{\ast}( \udn{T}; \mathbb{Z})$ is the exterior face ring of a  
simplicial complex (as implicitly described above). 
\end{conjecture}

Note that $I'$ is in fact quite different from the ideal $I$ in the statement of Theorem \ref{biggie}.
The conjecture is motivated by the fact that
$$ \Lambda[K] / I \cong \Lambda[K] / I' $$
in all of the existing computations.  These computations become very long even for small trees $T$ and integers $n$,
however, so the conjecture is ultimately based on only a few examples.

Robert Ghrist \cite{Ghr}  asked 
if all braid groups of planar graphs are right-angled Artin groups.  (Strictly speaking, he conjectured that all \emph{pure}
braid groups of planar graphs are right-angled Artin groups.  Here we follow our practice from \cite{FS2} and refer to either version
of the conjecture as Ghrist's conjecture.) A \emph{right-angled Artin group} is a group $G$ defined by
a presentation of the form $\langle x_1 , x_2 , \ldots, x_n \mid [ x_i , x_j ] \, \, (i,j) \in \mathcal{I} \rangle$, where
$\mathcal{I} \subseteq \{ 1, \ldots, n \} \times \{ 1, \ldots, n \}$ is an arbitrary indexing set.  
Equivalently, a right-angled Artin group $G_{\Gamma}$
may be defined by a simplicial graph $\Gamma$.  The generators of $G_{\Gamma}$ are vertices of $\Gamma$, and the generators $v_i$, $v_j$ commute
if they are adjacent in $\Gamma$.  
The class of right-angled
Artin groups thus includes, at opposite extremes, free groups (where $\mathcal{I} = \emptyset$ or $\Gamma$ is a discrete set) and 
free abelian groups (where $\mathcal{I} = \{ 1, \ldots, n \} \times \{ 1, \ldots, n \}$ or $\Gamma$ is a complete graph).  

There are well-known classifying spaces for right-angled Artin groups \cite{CD}.  Fix a right-angled Artin group $G_{\Gamma}$.  
Begin with a
copy of $S^{1}$, denoted $S^{1}_{i}$, having a single $0$-cell and a single $1$-cell, for each generator $v_i$.
Form the product $\prod_{i=1}^{n} S^{1}_{i}$.  Each open $j$-cell $e$ in the product corresponds naturally to a collection
of $j$ generators, one for each $1$-cell which contributes a factor to the $j$-cell $e$.  We throw out an open $j$-cell if and only if 
the corresponding set of generators fails to form a copy of the complete graph on $j$ vertices in $\Gamma$.  The resulting complex
$K_{\Gamma}$ is a $K(G_{\Gamma}, 1)$ complex.  It is easy to describe the cohomology ring of $K_{\Gamma}$ (and, thus, the cohomology
ring of $G_{\Gamma}$) using this complex (see \cite{Hatcher}):  $H^{\ast}( K_{\Gamma}; \mathbb{Z})$ is the 
exterior face ring of the largest simplicial
complex $K$ having $\Gamma$ as its $1$-skeleton.  A simplicial complex that is determined by its $1$-skeleton in this sense is called \emph{flag}.
      
It follows that Conjecture \ref{newconj} is a weakened version of Ghrist's conjecture.  Indeed, the latter conjecture (which was disproved in
\cite{FS2}) would have implied (as above)
that $H^{\ast}( \udn{\Gamma}; \mathbb{Z})$ is the exterior face ring of a flag complex.  A version of 
Conjecture \ref{newconj}
is also open for planar graphs, i.e., it is conceivable that the cohomology ring of $\udn{\Gamma}$ is an exterior face algebra for $\Gamma$ 
planar and $n$ arbitrary, although the current positive evidence is not strong.  We note finally that  
the cohomology rings of non-planar graph braid groups are not necessarily exterior face rings, since $\ud{2}{K_{5}}$ and $\ud{2}{K_{3,3}}$
are both homotopy equivalent to non-orientable surfaces \cite{A}. 


\section{Examples} \label{exples}


We conclude with some sample calculations of cohomology rings
$H^{\ast} \left( \ud{n}{T}; \mathbb{Z} \right)$.  It will be helpful to have a few 
lemmas for simplifying the calculations.

\begin{lemma} \label{dir}
Let $\phi$ be a standard cocycle, let $e$ be a bad edge of $\phi$, and let
$v = \iota(e)$.  For $i=0,1, \ldots, d(v)-1$, let $C_{i}$ denote the cloud
adjacent to $v$ in the direction $i$.  Suppose that the edge $e$ points in
direction $k>0$ from $v$.  
If $f( C_i ) = 0$ for all $i \neq k$, then $\phi$ is cohomologous to $0$.
\end{lemma}
\begin{proof}
The description of $\delta( R_{\mathcal{C}, v})$ implies that $\phi$ is
the only standard cocycle occurring in the sum
$\Sigma_{\mathcal{C},v, \tau} - \Sigma_{\mathcal{C}, v, \iota}$.
\end{proof}

\begin{lemma} \label{count}
Let $\phi$, $\psi$ be standard cocycles on $\ud{n}{T}$.  If
$[ \phi ] \cup [\psi ] \neq 0$ (where the brackets denote cohomology classes), 
then for any cloud $C$ of $\phi$,
$$ f_{\phi}(C) \geq \left[ \sum_{C' \subseteq C} f_{\psi}(C') \right]
+ | \{ e \in E( \psi ) \mid e \subseteq C \} |,$$
where the sum on the right runs over all clouds $C'$ of $\psi$ that are contained in $C$.
\end{lemma}

\begin{proof}
If $[ \phi ] \cup [\psi] \neq \emptyset$, then there is some standard cocycle 
$\nu$ such that each of the standard cocycles $\phi, \psi$ may be obtained from $\nu$
by breaking certain collections of edges in the cloud picture for $\nu$ and combining clouds.  The fact that
the cup product $[\phi] \cup [\psi]$ is non-zero also implies that
$E(\phi) \cap E(\psi) = \emptyset$. 

The cloud picture for $\phi$ is obtained from $\nu$ by breaking each edge $e \in E( \psi ) \subseteq E( \nu )$
and combining all clouds of $\nu$ incident with these edges.  The description of $ \nu $ implies that
each cloud of $ \nu $ is contained in one of the clouds of $ \phi $.  If we fix a cloud $C$ of $ \phi $ 
and break all edges
$e \in E( \psi) \subseteq E( \nu )$, then after combining all of the clouds $C'$ of $ \nu $ that are contained in $C$ 
there are precisely
$$f_{ \phi }(C) = \left[ \sum_{C' \subseteq C} f_{ \nu } (C') \right] + | \{ e \in E(\psi) \mid e \subseteq C \} |,$$
vertices in the cloud $C$,  
where the sum on the right runs over all clouds $C'$ of $ \nu $ that are contained in $C$.

We consider only the clouds $C'$ of $ \nu $ that are also clouds of $ \psi $.  Omitting the others (if any) from the above equality, 
we get
$$ f_{ \phi }(C) \geq \left[ \sum_{C' \subseteq C} f_{ \nu }(C') \right] + | \{ e \in E(\psi) \mid e \subseteq C \}|,$$
where the sum runs over all clouds $C' \subseteq C$ which are clouds of both $ \psi $ and $ \nu $.

We claim that $f_{ \nu }(C') = f_{ \psi }(C')$ for all such clouds $C'$; the lemma clearly follows from this claim.  
The statement that $C'$ is a cloud of
$ \nu $ and $ \psi $ implies that $C'$ touches only edges of $ \psi $.   
Now the required equality follows from the fact that $ \psi $ is obtained from $ \nu $ by breaking edges of $ \phi $ and combining
the clouds; this procedure doesn't change the number of vertices in $C'$ by our assumptions, thus proving the claim.
\end{proof}

\begin{lemma} \label{principle}
Let $\phi$ be a critical $1$-cocycle defined on $\ud{n}{T}$; and let
$[ \phi ]$ denote its cohomology class.  If $[ \phi ] \cup [ \theta ]  \neq 0$ for
some other cohomology class $[ \theta ] \in H^{j} \left( \udn{T} \right)$ ($j \geq 1$), 
then there is some cloud $C$ of $\phi$ such that:
\begin{enumerate}
\item $C$ contains at least one essential vertex of $T$, and
\item $f_{\phi}(C) \geq 2$.
\end{enumerate}
\end{lemma}

\begin{proof}
Fix a critical $1$-cocycle $\phi$.  If there is a cohomology class $[ \theta ]$ as
in the statement of the Lemma, then there is some critical $1$-cocycle $\theta'$ such 
that $[ \phi ] \cup [ \theta' ] \neq 0$, since $H^{\ast}(\udn{T})$ is generated by 
$H^{1}(\udn{T})$, and a basis for $H^{1} (\udn{T})$ is determined by the cohomology classes
of critical $1$-cocycles, by Corollary \ref{surjective}.         

Since $[ \phi ] \cup [ \theta' ] \neq 0$, there is some $[ \psi ] \in H^{2}(\udn{T}) - \{ 0 \}$
where $\psi$ is a standard $2$-cocycle (which need not be critical) and the equivalence class $[\psi]$ is the least
upper bound of $\{ [ \phi ] , [ \theta' ] \}$.  Let $e$ be the edge
of $\psi$ which is not an edge of $\phi$.  Now $[ \psi ] \neq 0$, so, by Lemma \ref{dir}, at least one
of the clouds $C_1$ incident with $e$ satisfies $f_{\psi}(C_{1}) > 0$.  It follows that $\phi$, which is
the result of breaking the edge $e$ in $\psi$ and combining clouds, satisfies $f_{\phi}(C) \geq 2$, where 
$C$ is the unique cloud of $\phi$ containing $e$.  This proves property (2).  We have that $e$ is the 
edge of $\theta'$, and since this is a critical cocycle, $\iota(e)$ must be essential.  It follows that
$C$ also satisfies (1).
\end{proof}

\begin{example} \label{modh2a} 
We consider $H^{\ast}( \ud{4}{T} )$, where $T$ is the tree depicted Figure \ref{modh},
with the given choices of embedding and basepoint $\ast$.

\begin{figure} [!h]
\begin{center}
\includegraphics{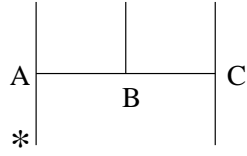}
\caption{This is the tree $T$ from Example \ref{modh2a}.  The essential vertices are
labelled $A$, $B$, and $C$.}
\label{modh}
\end{center}
\end{figure}

A critical $1$-cocycle $\phi$ on $\ud{4}{T}$ is a standard cocycle
satisfying two conditions:  i)  the (unique) edge $e$ of $\phi$ is one
of $e_{A}$, $e_{B}$, $e_{C}$, where $e_{X}$ is the unique edge such that
$\iota( e_{X} ) = X$ and $e_{X}$ points in direction $2$ from $X$; ii) the
cloud $C$ in direction $1$ from $\iota(e)$ satisfies $f(C) > 0$.  Theorem
\ref{biggie} implies that the critical $1$-cocycles determine a distinguished basis
for $H^{1}( \ud{4}{T} )$.

Fix a particular critical $1$-cocycle $\phi$ and let $[ \phi]$ denote its
cohomology class.  If $[\phi]$ cups non-trivially with some other member $[ \psi ]$
of the distinguished basis (where $\psi$ is a critical $1$-cocycle),
then it follows from Lemma \ref{principle} that some cloud $C$ of $\phi$ contains
at least one essential vertex and satisfies $f_{\phi}( C ) \geq 2$.  
It follows that there
are exactly four critical $1$-cocycles $\phi$ having non-trivial cup products; these cocycles label
the vertices in Figure \ref{modh2}.

\begin{figure} [!h]
\begin{center}
\includegraphics{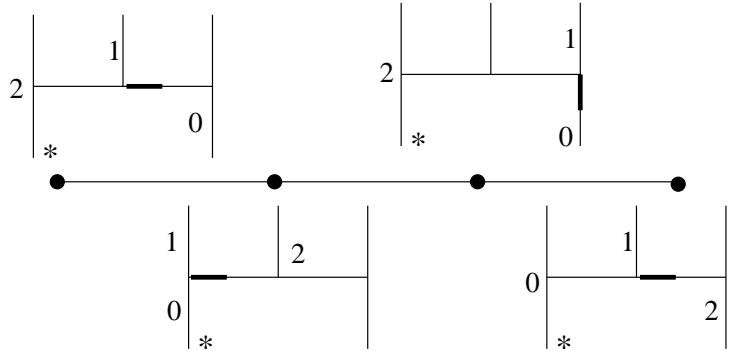}
\caption{Here is a portion of a certain simplicial complex $K$. 
The remaining parts of $K$ are isolated vertices.  The cohomology ring $H^{\ast}(\udn{T})$ is the exterior
face ring defined by $K$.}  
\label{modh2}
\end{center}
\end{figure}

The edges are to be labelled by ($2$-cocycle representatives of) the cup
products of their endpoints.  The precise details of this labelling are an 
exercise.  The cup product of any other pair of elements
of the standard basis for $H^{1}( \ud{4}{T})$ is $0$.  It follows from
Proposition 4.1 of \cite{Me} that
$H^{1}( \ud{4}{T} ) \cong \mathbb{Z}^{18}$ and
$H^{2}( \ud{4}{T} ) \cong \mathbb{Z}^{3}$ as abelian groups.

It now follows that $H^{\ast}(\ud{4}{T} ; \mathbb{Z}) \cong \Lambda(K)$,
where $K$ is the simplicial complex consisting of a line segment that is
divided into $3$ pieces (as in Figure \ref{modh2}) and $14$ isolated vertices. \qed
\end{example} 

\begin{example}
Consider now the cohomology ring $H^{\ast}( \ud{5}{H} )$, where $H$ is the 
graph that is homeomorphic to the letter ``H".  Proposition 4.1 from \cite{Me} 
implies that
$H^{1}( \ud{5}{H} ; \mathbb{Z} ) \cong \mathbb{Z}^{20}$,
$H^{2}( \ud{5}{H} ; \mathbb{Z} ) \cong \mathbb{Z}^{5}$, and 
$H^{n}( \ud{5}{H} ; \mathbb{Z} ) \cong 0$ for $n>2$.  We let $A$ denote
the left-hand essential vertex and $B$ denote the right-hand essential vertex.
Our basepoint $\ast$ is the vertex of degree $1$ at the bottom left.

\begin{figure} [!h]
\begin{center}
\includegraphics{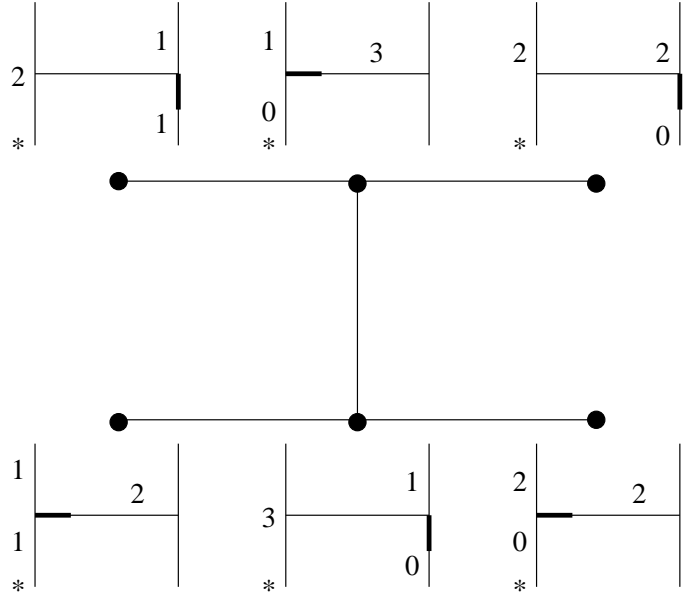}
\caption{This is part of the simplicial complex $K$.  The rest of $K$ consists of isolated vertices.
The cohomology ring $H^{\ast}(\ud{5}{H})$ is the exterior face ring determined by $K$.
}
\label{h}
\end{center}
\end{figure}

The critical $1$-cocycles $\phi$ of $\ud{5}{H}$ (which represent a
distinguished basis for $H^{1} ( \ud{5}{H} ; \mathbb{Z} )$) are the
standard $1$-cocycles determined by the following two properties:  i)
the (unique) edge $e$ of $\phi$ satisfies $\iota(e) = A$ or $\iota(e) = B$
and points in direction $2$ from $\iota(e)$;  ii) the cloud $C$ in direction
$1$ from $\iota(e)$ satisfies $f_{\phi}(C) > 0$.

If $[ \phi ]$ cups non-trivially with some cohomology class $[ \psi ]$
(where $\phi$, $\psi$ are both critical $1$-cocycles), then, by Lemma \ref{principle}, 
there must be some cloud $C'$ of $\phi$ which contains
some essential vertex and such that $f_{\phi}(C') \geq 2$.  This leaves only
six elements of the distinguished basis for $H^{1} ( \ud{5}{H})$ which
have non-trivial cup products; cocycle representatives for these cohomology
classes label the vertices in Figure \ref{h}.

The edges should be labelled by ( cocycle representatives of) the corresponding
cup products.  Note that the edge labels represent a basis for
$H^{2}( \ud{5}{H} )$.  We conclude that $H^{\ast}( \ud{5}{H} ) \cong
\Lambda (K)$, where $K$ is the simplicial complex consisting of $14$
isolated vertices and the letter ``H" (as cellulated in Figure \ref{h}). \qed
\end{example}

\begin{example} \label{mess}
Let $T$ be the tree depicted in Figure \ref{plus}.  It is a consequence of
Theorem 3.7 from \cite{Me} that $H^{0}(\ud{4}{T}) \cong 0$, $H^{1} (\ud{4}{T}) \cong \mathbb{Z}^{50}$,
$H^{2} (\ud{4}{T}) \cong \mathbb{Z}^{18}$, and $H^{n}(\ud{4}{T}) \cong 0$ for $n >2$.

\begin{figure}[!h]
\begin{center}
\includegraphics{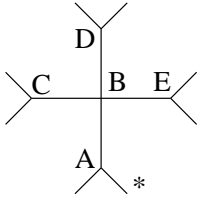}
\caption{Here is the tree $T$ that is under consideration in Example \ref{mess}.}
\label{plus}
\end{center}
\end{figure}

We compute the cohomology ring 
$H^{\ast} \left( \ud{4}{T} ; \mathbb{Z} \right)$.  We will need to consider
a large number of cocycles in our argument, so it will be useful to 
introduce new notation.  We describe a critical $1$-cocycle $\phi$ by a
$4$- (or $5$-) tuple as follows.
\begin{enumerate}
\item The first entry is $\iota(e)$, where $e$ is the unique edge of
the cocycle $\phi$.  By the description of critical cocycles, it follows
that $\iota(e) \in \{ A, B, C, D, E \}$. 
\item The last $3$ (or $4$) entries are the numbers of cells lying in,
respectively, the directions $0$, $1$, $2$, and (possibly) $3$ from
$\iota(e)$.  We therefore have a $4$-tuple if $\iota(e) \in \{ A, C, D,
E \}$; $5$-tuple if $\iota(e) \in \{ B \}$.  Note that this cell count
includes the edge $e$.
\item Finally, we indicate the location of the edge $e$ by placing a bar
over the $(i +2)$nd entry of the tuple, if the direction from $\iota(e)$
to $\tau(e)$ is $i$.
\end{enumerate}
For instance, if we use this notation to describe the six critical cocycles 
which label the vertices in Figure \ref{h}, then, reading from the top
left, we get
$$ (B, 2, 1, \bar{2}); (A, 0,1,\bar{4}); (B,2,2,\bar{1});
(A,1,1,\bar{3}); (B,3,1,\bar{1}); (A,0,2,\bar{3}).$$
We use a somewhat modified notation in Figure \ref{big}.  For instance,
the $4$-tuple $(A, 0, 1, \bar{3})$ is replaced with a circled ``A"; the
second, third, and fourth entries of the $4$-tuple are written in clockwise
order around the interior of the circle, beginning from below the ``A".
The ``3" is circled to indicate the location of the edge.

A vertex represents the cohomology class determined by its label.  Two
vertices are connected by an edge if they represent $1$-cocycles having
a common upper bound.  In almost every case, the upper bound in question
is a critical cocycle.  There are exactly $7$ edges representing non-critical
cocycles (these edges are 
dotted in the picture; the numbering of the edges is keyed to the enumerated points below). 
We consider each case in turn.  Let, for instance,
$(A, 0, 1, \bar{3})-(E,2,1,\bar{1})$ denote the edge running between the
indicated vertices.

\begin{figure} [!h]
\begin{center}
\includegraphics{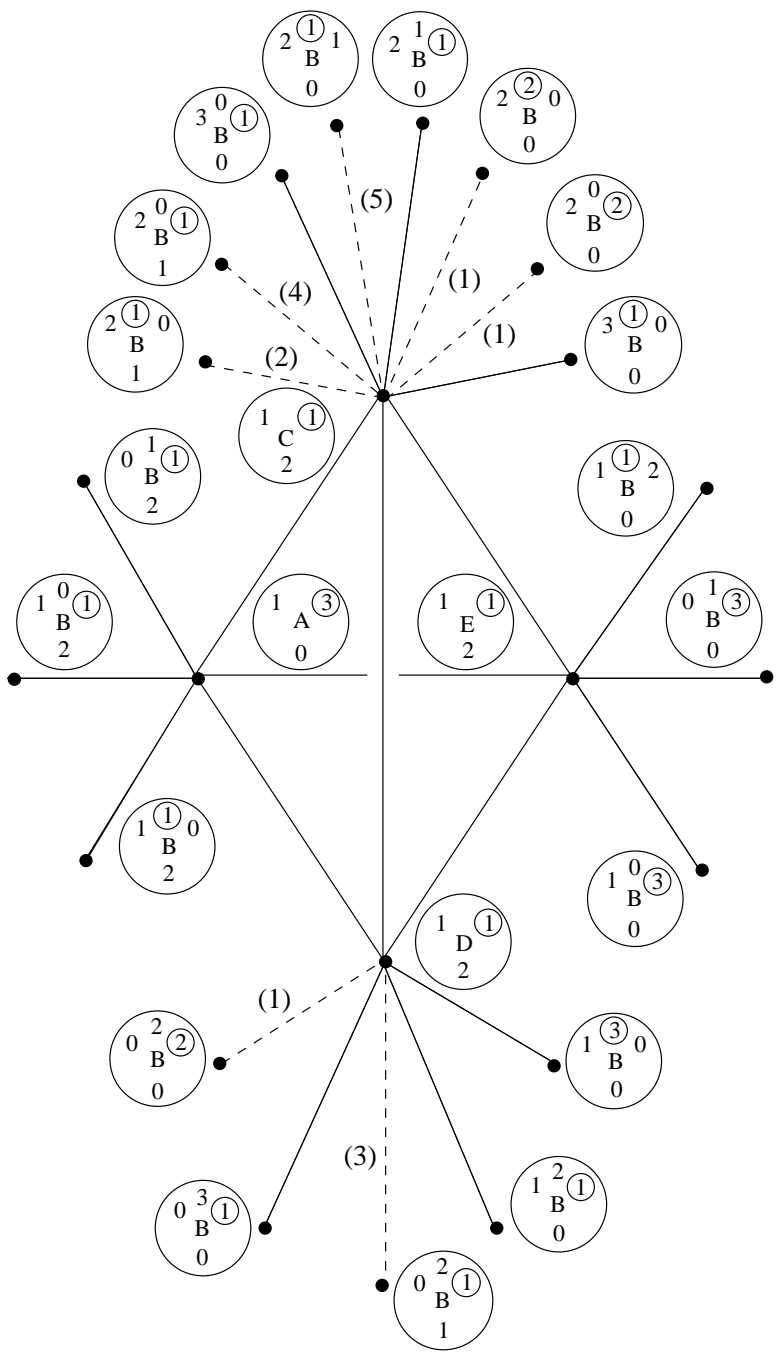}
\caption{Here is part of the simplicial complex $K$ which determines the cohomology ring $H^{\ast}(\ud{4}{T})$.  The seven dotted
edges can be eliminated by a change of basis in $H^{\ast}(\ud{4}{T})$.}
\label{big}
\end{center}
\end{figure}

\begin{enumerate}
\item \underline{$(B,0,0,2,\bar{2})-(D,2,1,\bar{1})$, 
$(B,0,2,0,\bar{2})-(C,2,1,\bar{1})$, $(B,0,2,\bar{2},0)-$} \newline \underline{$(C,2,1,\bar{1})$}
Each of these edges represents a standard $2$-cocycle which satisfies the
hypotheses of Lemma \ref{dir}.  Accordingly, all of the cup products in question
are $0$.
\item \underline{ $(B,1,2,\bar{1},0)-(C,2,1,\bar{1})$}  The standard 
$2$-cocycle represented by the given edge has the form in Figure \ref{ex2}.

\begin{figure}[!h]
\begin{center}
\includegraphics{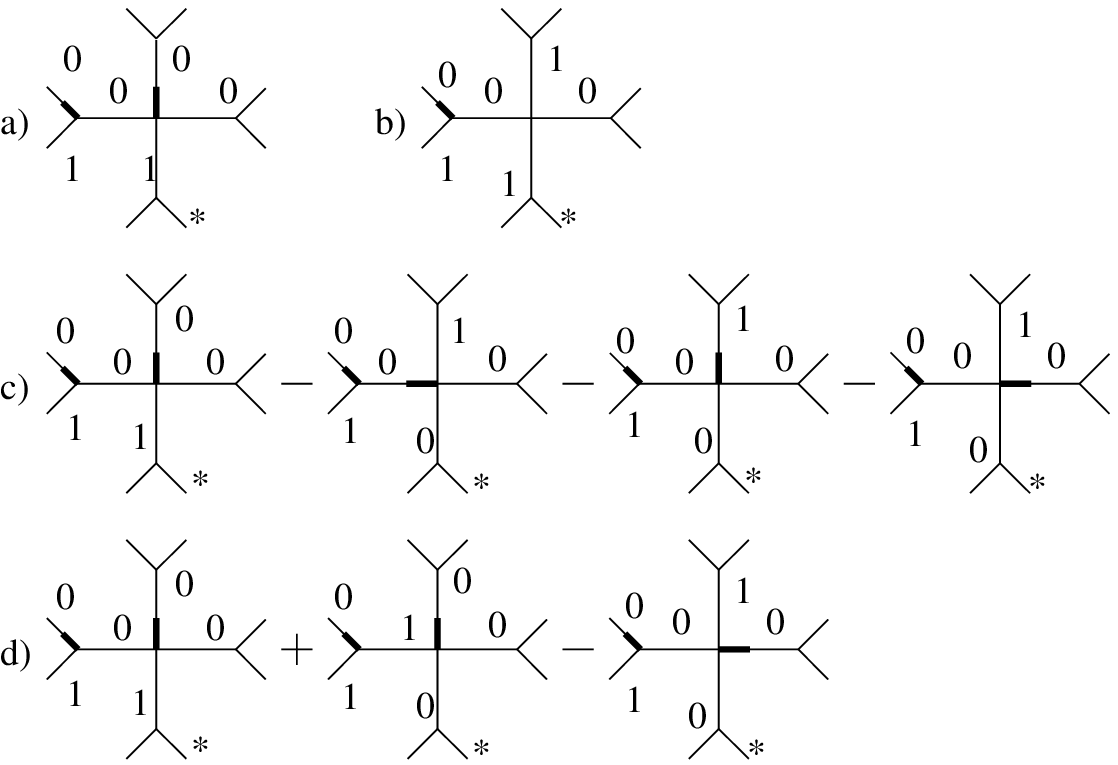}
\caption{a) This standard cocycle $\phi$ is the result of cupping 
$(B,1,2,\bar{1}, 0)$ and $(C,2,1,\bar{1})$. b) This is the cochain $R_{\mathcal{C},B}$.
c) This is the coboundary $\delta( R_{\mathcal{C},B})$. d)  This is the result of 
applying the procedure from Theorem \ref{biggie} to the final three terms of the sum from c).}
\label{ex2}
\end{center}
\end{figure}

The bad edge $e$ of the standard cocycle $\phi$ with cloud picture 
$\mathcal{C}$ depicted in Figure \ref{ex2} a) satisfies $\iota(e) = B$.  Figure
\ref{ex2} b) depicts the cochain $R_{\mathcal{C},B}$, where the central vertex of
degree $4$ belongs to the same cloud as the basepoint $\ast$.  In 
Figure \ref{ex2} c), we see $\delta( R_{\mathcal{C},B})$.  If we apply the
rewrite procedure from Theorem \ref{biggie} to the final $3$ terms in c) (and leave
the first term alone), we arrive at the sum of standard cocycles in d).
(The precise details in this last step are an exercise.)  The sum $\Sigma$
in d) is clearly a sum of coboundaries -- this is a direct consequence of the way we computed it --   
so $\Sigma = 0$ on cohomology.  It
follows directly that the cup product of the cohomology classes
$(B,1,2,\bar{1},0)$ and $(C,2,1,\bar{1})$ is represented by $\phi - \Sigma$,
and $\phi - \Sigma$ is a sum of critical cocycles.

We would like to choose a basis for $H^{\ast}(\ud{4}{T}; \mathbb{Z})$ that
will give the simplest possible formulas for the cup product.  To achieve
this, consider the result of breaking the left-hand edge in each of the
cloud pictures of d).  We arrive at the sum 
$\Sigma' = (B,1,2,\bar{1},0) + (B,3,\bar{1},0) - (B,2,1,\bar{1})$ of
critical cocycles.  It is straightforward to check that the cup product
of $[ \Sigma' ]$ with any other cohomology class of dimension at least $1$ is $0$.  Indeed, 
by Lemma \ref{count}, it suffices to check that the cup product of $[ \Sigma' ]$
with $[ (C,2,1,\bar{1})]$ is $0$, but the product in question is represented
by the sum $\Sigma$ in Figure \ref{ex2} d), which is $0$ on cohomology.  We replace
the class $[(B,1,2,\bar{1},0)]$ with 
$[ (B,1,2,\bar{1},0) + (B,0,3,\bar{1},0) - (B,0,2,1,\bar{1})]$ in our
standard basis.
\item \underline{ $(B,1,0,2,\bar{1})-(D,2,1,\bar{1})$}  This edge is labelled
by the following cloud picture:

\begin{figure}[!h]
\begin{center}
\includegraphics{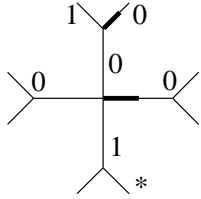}
\caption{This is the standard cocycle $\phi$ representing the cup product of $(B,1,0,2, \bar{1})$ with $(D,2,1,\bar{1})$.}
\label{final}
\end{center}
\end{figure}

This case follows the pattern of 2); we apply the rewrite procedure from
Theorem \ref{biggie} to express the above cocycle $\phi$ as a sum of critical
cocycles $\Sigma''$, i.e., 
$$ \phi = \Sigma'',$$
where the equality is valid on the level of cohomology.  The sum 
$\phi - \Sigma''$ is sum of standard cocycles, each of which has a given
fixed edge $e$ such that $\iota(e) = D$.  If we break this edge in each
standard cocycle of the sum $\phi - \Sigma''$ and combine clouds, then we
arrive at the sum $\widehat{\Sigma} = 
(B,1,0,2,\bar{1}) + (B,0,1,2,\bar{1}) + (B,0,0,3,\bar{1})$.  This represents
a cohomology class which, as in 2), has no non-zero cup products.  We
replace $[(B,1,0,2,\bar{1})]$ with $[ \widehat{\Sigma} ]$ in the standard
basis.
\item \underline{ $(B,1,2,0,\bar{1})-(C,2,1,\bar{1})$}  
This case again follows
the pattern of 2) and 3).  We replace the cohomology class 
$[ (B,1,2,0,\bar{1})]$ in the standard basis with
$[ (B,1,2,0,\bar{1}) + (B,0,3,0,\bar{1}) + (B,2,1,\bar{1})]$.
\item \underline{ $(B,0,2,\bar{1},1)-(C,2,1,\bar{1})$}  Replace
$[(B,0,2,\bar{1},1)]$ with $[ (B,0,2,\bar{1},1) + (B,0,2,1,\bar{1})]$.
\end{enumerate}
 
Items 1)-5) show that we can arrange, after a change of basis, for the seven
dashed edges in Figure \ref{big} to be deleted.  It is routine to check that 
the $2$-dimensional cohomology classes which should label the edges in the remaining graph
form a basis for $H^{2}( \ud{4}{T})$.  
The cup product ring is thus the
exterior face ring of the graph consisting of a copy of $K_4$ with three additional edges 
attached at each vertex, along with $34$ isolated vertices.     
\qed

\end{example}
 

\bibliography{bib}
\nocite{*}
\bibliographystyle{plain}

\end{document}